\documentclass{article}
\usepackage{authblk}
\usepackage[utf8]{inputenc}
\usepackage{amsmath}
\usepackage{amsbsy}
\usepackage{amssymb}
\usepackage{graphicx}
\usepackage{dsfont}
\usepackage{upgreek}
\usepackage{textcomp}
\usepackage{braket}
\usepackage[margin=1in]{geometry}
\usepackage{amsthm}
\usepackage{mathrsfs}
\usepackage{mathtools}
\usepackage[table,svgnames]{xcolor}
\usepackage{graphicx}
\usepackage{tikz}
\usepackage{float}
\usepackage{enumerate}
\usetikzlibrary{arrows}
\usepackage[toc,page]{appendix}
\usepackage{cleveref}
\usepackage[toc,page]{appendix}
\crefname{appsec}{appendix}{appendices}
\usepackage{xcolor}

\crefname{appsec}{appendix}{appendices}
\numberwithin{equation}{section}

\newtheorem{theorem}{Theorem}[section]
\newtheorem{lemma}{Lemma}[section]

\newtheorem{remark}{Remark}[section]

\newtheorem{proposition}{Proposition}[section]

\providecommand{\keywords}[1]
{
  \textbf{Key words.} #1
}
\providecommand{\classification}[1]
{
  \textbf{Mathematics Subject Classification.} #1
}

\title{Remark on the Adiabatic Limit of Quantum Zakharov System}
\author{Brian Choi\thanks{Email: \texttt{choigh@bu.edu}}%
  }
\affil{Department of Mathematics and Statistics, Boston University, Boston, MA 02215, USA}
\date{}
\begin{document}
\maketitle\vspace{-7ex}
\begin{abstract}
This paper is concerned with the low regularity well-posedness of the adiabatic limit of the quantum Zakharov system, which is a modified nonlinear Schr\"odinger equation (NLSE). As the quantum parameter tends to zero, the modified NLSE formally converges to the standard NLSE, and we show this rigorously.
\end{abstract}
\\
\classification{35B30,35Q40,35Q55}
\label{}

\section{Introduction} 
Consider \begin{equation}\label{adiabatic}
\begin{split}
\begin{cases}
    i\partial_t u + \Delta u -\epsilon^2\Delta^2u=-(I-\epsilon^2\Delta)^{-1}(|u|^2)u,\: (x,t)\in M\times \mathbb{R}\\
    u(0)=u_0\in H^s(M),
    \end{cases}    
\end{split}
\end{equation}
for $M = \mathbb{R}$ or $\mathbb{T}$ where $J_\epsilon:=(I-\epsilon^2\Delta)^{-1}$. Smooth solutions to \eqref{adiabatic} have two notable conserved quantities:
\begin{equation}\label{conservation}
\begin{split}
    \text{MASS}&=\lVert u(t)\rVert_{L^2}^2,\\
    ENERGY&=E[u(t)] = \frac{\epsilon^2}{2}\lVert \partial_{xx}u\rVert_{L^2}^2 + \frac{1}{2}\lVert \partial_x u \rVert_{L^2}^2 -\frac{1}{4}\lVert J_\epsilon(|u|^2)|u|^2\rVert_{L^1}.
\end{split}
\end{equation}
This analysis of \eqref{adiabatic} is motivated by the quantum Zakharov system 
\begin{equation}\label{qzs}
    \begin{cases}
    i\partial_t E + \Delta E-\epsilon^2\Delta^2E=nE,\: (x,t)\in \mathbb{R}^d\times \mathbb{R}\\
    \frac{1}{\lambda^2}\partial_{tt}n-\Delta n+\epsilon^2\Delta^2n=\Delta(|E|^2),\\
    E(0)=E_0,\: n(0)=n_0,\: \partial_t n(0)=n_1.
    \end{cases}
\end{equation}
The classical Zakharov system \cite{zakharov1972collapse} (\eqref{qzs} with $\epsilon=0$) describes the propagation of Langmuir waves in an ionized plasma where the complex-valued $E(x,t)$ describes a slowly-varying envelope of a rapidly oscillating electric field and the real-valued $n(x,t)$ describes the deviation of the ion density from its mean. On the other hand, \eqref{qzs} accounts for the \textit{quantum effects} on the nonlinear interaction between $E$ and $n$ with $\epsilon = \frac{\hbar w_i}{k_B T_e}$ where the fourth-order modification is experimentally relevant in the presence of dense and cold plasmas, which occurs in astrophysical settings; see \cite{garcia2005modified,haas2011quantum,haas2009quantum} for more physical background.  

The well-posedness theory of the classical Zakharov system is a mature subject by now where \cite{ozawa1992existence,kenig1995zakharov,bourgain1996wellposedness,ginibre1997cauchy,colliander2008low,bejenaru20092d,bejenaru2011convolutions} is a glimpse of references on the topic. On the other hand, the study of \eqref{qzs} is relatively recent: \cite{guo2013global,1078-0947_2016_10_5445,fang2016fourth,chen2017low,fang2019semi}. A central theme in this wave of recent work is the convergence of \eqref{qzs} as $\epsilon\rightarrow 0$ (semi-classical limit) or $\lambda \rightarrow \infty$ (adiabatic limit). The latter, in the context of classical Zakharov system, was studied by \cite{schochet1986nonlinear,ozawa1992nonlinear} where the limiting dynamics was rigorously identified, with an optimal convergence rate, as the focusing cubic nonlinear Schr\"odinger equation (NLSE). Similarly when $\epsilon>0$, \cite{fang2016fourth} showed the adiabatic limit of \eqref{qzs} to \eqref{adiabatic}. Our work is concerned with solving \eqref{adiabatic} with infinite-energy data, thereby extending \cite[Proposition 2.5]{fang2016fourth} that used \eqref{conservation} to obtain coercive bounds in $H^2$. Our proof adopts the method of Fourier restricted norms, first initiated by Bourgain \cite{bourgain1993fourier}. It was also used to obtain the low regularity well-posedness of KdV and quadratic NLSE on $\mathbb{R}$ in \cite{kenig1996bilinear,kenig1996quadratic}, respectively.

\begin{theorem}\label{thm311}
For all $s\geq 0$, \eqref{adiabatic} is globally well-posed in $H^s(M)$. Furthermore the data-to-solution map is locally Lipschitz in $H^s(M)$.
\end{theorem} 

For $s<0$, it is shown that the solution map for \eqref{adiabatic} on $\mathbb{T}$, if it exists, fails to be uniformly continuous by constructing a family of explicit examples in the spirit of \cite{burq2002instability}. Consequently $s=0$ is the threshold above which the contraction mapping argument yields well-posed solutions. Topics such as the ill-posedness on $\mathbb{R}$ and the blow-up dynamics of \eqref{adiabatic} as $\epsilon\rightarrow 0$ are not covered in this work. 

Having established the existence of global solutions, the semi-classical limit of \eqref{adiabatic} to NLSE is discussed. For this limit, the growth of Sobolev norms of solutions needs to be uniform in $\epsilon$; in particular, the smoothing estimates developed for \eqref{thm311} are to no avail. It is shown that the solution to \eqref{adiabatic} converges to the solution to the cubic NLSE on $[-T,T]$ for every $T>0$ whereas the convergence fails for $T=\infty$.

\begin{theorem}\label{thmII}
Let $u^{(\epsilon)}$ and $u$ be the solutions to \eqref{adiabatic} and cubic NLSE with data $u_0^{(\epsilon)}$ and $u_0$, respectively. Let $\{u_0^{(\epsilon)}\}_\epsilon$ be bounded in $H^2(M)$ and $u_0^{(\epsilon)}\xrightarrow[\epsilon\rightarrow 0]{} u_0$ in $H^1(M)$. Then $u^{(\epsilon)}\xrightarrow[\epsilon\rightarrow 0]{}u$ in $C([-T,T];H^1(M))$ for every $T>0$.
\end{theorem}

We outline the organization of this paper. In \cref{section wp}, \eqref{thm311} is proved separately for $\mathbb{R},\mathbb{T}$. On $\mathbb{R}$, the smoothing estimate \eqref{mainprop} is shown by direct estimation using the Cauchy-Schwartz inequality. On $\mathbb{T}$, an appropriate bound in $L^4(\mathbb{T}\times\mathbb{R})$ is obtained by a counting argument. In \cref{section limit}, the limit of \eqref{adiabatic} as $\epsilon\rightarrow 0$ is considered. Explicit examples are given to show that the uniform convergence does not hold globally in time.

\section{Well-posedness}\label{section wp}

Let $w(\xi) = \epsilon^2 \xi^4 + \xi^2$ where $\epsilon \geq 0$. Denote $\xi \in \mathbb{R},\:k \in \mathbb{Z},\: \tau \in \mathbb{R}$ as the dual variables of $x\in\mathbb{R},\: x\in \mathbb{T},\: t \in \mathbb{R}$, respectively where the notation $\widehat{u}$ is used for either space Fourier transform/series or spacetime Fourier transform/series, depending on the context. Let $\langle \xi \rangle = (1+|\xi|^2)^{1/2}$. For $s,b \in \mathbb{R}$, define
\begin{equation*}
    \lVert u \rVert_{X^{s,b}}= \lVert \langle \xi \rangle^s \langle \tau + w(\xi)\rangle^b \hat{u}(\xi,\tau)\rVert_{L^2_{\xi, \tau}},
\end{equation*}
and its restricted norm $\| u \|_{X^{s,b}_T} = \inf\limits_{\Tilde{u} = u,\: t \in [0,T]} \| \Tilde{u}\|_{X^{s,b}}$ where by the time-reversal symmetry, it suffices to solve \eqref{adiabatic} on $[0,T]$. The nonlinear term $N(u) = J_\epsilon(|u|^2)u$ undergoes smoothing under our choice of norm. The proof of \cref{mainprop}, given in the Appendix, requires a technical analysis of a cubic function in the Fourier space.
\begin{proposition}\label{mainprop}
Let $s \in [0,\infty), \gamma\in [\frac{1}{3},\frac{1}{2}),T \in (0,1]$ and $a\in [0,\frac{4}{3})$. Then for every $\epsilon>0$ and $b\in (\frac{1}{2},1-\gamma)$,
\begin{equation*}
\lVert J_\epsilon(u\overline{v})w \rVert_{X^{s+a,-\gamma}_{T}}\lesssim_{s,a,\gamma,b,\epsilon} \lVert u \rVert_{X^{s,b}_{T}}\lVert v \rVert_{X^{s,b}_{T}}\lVert w \rVert_{X^{s,b}_{T}}.    
\end{equation*}
\end{proposition}
\begin{remark}
The proof of \cref{mainprop} is adopted from \cite[Lemma 2.4]{kenig1996bilinear} where the analysis of the KdV equation (third-order derivative) leads to a direct estimation of a quadratic polynomial in the Fourier space; on the other hand, \eqref{adiabatic} (fourth-order derivative) demands an appropriate bound of a third-order polynomial. This process is technical and thus is postponed to the Appendix.
\end{remark}
\begin{remark}
Our trilinear estimates in the $X^{s,b}$ space defined on $\mathbb{R}$ or $\mathbb{T}$ have subtle differences. The derivative gain in \cref{mainprop} is in $s$ quantified by $a>0$. Since the rise in $b$ is $(1-\gamma)-(-\gamma)=1$, there is no derivative gain in $b$ by the fourth statement of \cref{l:32}. On the other hand, the derivative gain for \cref{lemmaembedding} is not in $s$ but in $b$, since the rise, as mentioned previously, is $\frac{5}{8}<1$.   
\end{remark}
The mapping properties of the following operators are repeatedly used in our analysis.
\begin{lemma}\cite[Section 1.5]{cazenave2003semilinear}\label{cazenave}
Let $s\in \mathbb{R}, \epsilon>0$ and $X=H^s(\mathbb{R}^d)\:\text{or}\:L^p(\mathbb{R}^d)$ for $p\in [1,\infty)$. Then $\lVert J_\epsilon f \rVert_{X}\leq \lVert f \rVert_{X}$ for all $f\in X$. Moreover $J_\epsilon: H^s(\mathbb{R}^d)\rightarrow H^{s+2}(\mathbb{R}^d)$ is bounded with the best constant $\leq \epsilon^{-2}$.
\end{lemma}
\begin{lemma}\label{l:32} Let $T \in (0,1],\:s,b\in\mathbb{R}$. Let $\eta$ be a smooth bump function on $\mathbb{R}$. Then
\begin{enumerate}
    \item \cite[Lemma 2.8]{tao2006nonlinear}: $\lVert \eta(t) e^{-it w(\frac{\nabla}{i})}f\rVert_{X^{s,b}}\lesssim_{b,\eta}\lVert f \rVert_{H^s(M)}$.
    \item \cite[Corollary 2.10]{tao2006nonlinear}: For every $b > \frac{1}{2}$, we have
    \begin{align*}
        X^{s,b}_T \hookrightarrow C([0,T];H^s(M)).
    \end{align*}
    \item \cite[Lemma 2.11]{tao2006nonlinear}: Let $-\frac{1}{2}<b^\prime\leq b<\frac{1}{2}$. Then,
    \begin{align*}
    \lVert \eta(\frac{t}{T})u\rVert_{X^{s,b^\prime}}\lesssim_{b^\prime,b,\eta}T^{b-b^\prime}\lVert u \rVert_{X^{s,b}}.    
    \end{align*}
    \item \cite[Proposition 2.12]{tao2006nonlinear}: Let $b>\frac{1}{2}$. Then,
    \begin{align*}
    \left|\left| \eta(t) \int_0^t e^{-i(t-t^\prime)w(\frac{\nabla}{i})} F(t^\prime)dt^\prime \right|\right|_{X^{s,b}} \lesssim  \lVert F \rVert_{X^{s,b-1}}.    
    \end{align*}
\end{enumerate}
\end{lemma}
\begin{proof}[Proof of \cref{thm311} for $M=\mathbb{R}$]\label{lwp}
Let $s,\gamma,\delta,a,\epsilon$ and $b$ be as above. For $u_0\in H^s$, define
\begin{equation}\label{functionspace}
X=\{u\in X^{s,b}_{T}: \lVert u\rVert_{X^{s,b}_{T}}\leq 2C \lVert g \rVert_{H^s}\}    
\end{equation}
where $C$ satisfies $\lVert U_\epsilon(t)f\rVert_{X^{s,b}_{T}}\leq C \| f \|_{H^s}$ for all $f\in H^s$. Define
\begin{equation}\label{contraction}
\Gamma u=U_\epsilon(t)u_0-i\int_0^t U_\epsilon(t-t^\prime)N(u)(t^\prime)dt^\prime.  
\end{equation}
Since
\begin{equation*}
\begin{split}
    \lVert \Gamma u \rVert_{X^{s,b}_{T}}&\lesssim \lVert u_0 \rVert_{H^s}+\lVert N(u)\rVert_{X^{s,b-1}_{T}}\lesssim \lVert u_0 \rVert_{H^s}+T^{1-(b+\gamma)}\lVert N(u)\rVert_{X^{s,-\gamma}_{T}}\\
    &\lesssim \lVert u_0 \rVert_{H^s}+T^{1-(b+\gamma)}\lVert u \rVert_{X^{s,b}_{T}}^3\lesssim \lVert u_0 \rVert_{H^s}+T^{1-(b+\gamma)}\lVert u_0 \rVert_{H^s}^3
\end{split}
\end{equation*}
by \cref{l:32}, it follows that $\lVert \Gamma u \rVert_{X^{s,b}_{T}}\leq C\lVert u_0 \rVert_{H^s}+C_1T^{1-(b+\gamma)}\lVert u_0\rVert_{H^s}^3$, and hence $\Gamma:X\rightarrow X$ by taking $T<(\frac{C}{C_1}\| u_0 \|^{-2})^{\frac{1}{1-(b+\gamma)}}$.
For all $u,v \in X$,
\begin{equation*}
\begin{split}
    \lVert \Gamma u-\Gamma v \rVert_{X^{s,b}_{T}} &=\left|\left| \int_0^t U_\epsilon(t-t^\prime)(N(u)(t^\prime)-N(v)(t^\prime))dt^\prime\right|\right|_{X^{s,b}_{T}}\\
    &\lesssim \lVert N(u)-N(v)\rVert_{X^{s,b-1}_{T}}\lesssim T^{1-(b+\gamma)}\lVert N(u)-N(v) \rVert_{X^{s,-\gamma}_{T}}\\
    &\lesssim T^{1-(b+\gamma)}\lVert u_0 \rVert_{H^s}^2 \lVert u-v \rVert_{X^{s,b}_{T}},   
\end{split}
\end{equation*}
where the last inequality follows from
\begin{equation*}
    N(u)-N(v)=J_\epsilon(|u|^2)(u-v)+J_\epsilon(u\cdot(\overline{u-v}))v+J_\epsilon((u-v)\overline{v})v.
\end{equation*}
By shrinking $T$ if necessary, $\Gamma$ is a contraction on $X$. The Lipschitz continuous dependence on data is proved similarly.

To extend the solutions globally in time, it is shown that the Sobolev norm of solutions grows at most polynomially. We claim that for every $s \geq 0,\: \epsilon>0$, there exists a non-decreasing function $C_s^{(\epsilon)}:[0,\infty)\rightarrow [0,\infty)$ such that
\begin{equation*}
\lVert u^{(\epsilon)}(t)\rVert_{H^s}\leq C_s^{(\epsilon)}(\lVert u_0 \rVert_{H^s})\langle t \rangle^{\frac{1}{2}\Big(3^{\frac{3s}{4}+1}-1\Big)},    
\end{equation*}
for all $t \geq 0$. 

The proof is by induction on $k \geq 0$. Let $I_k = [\frac{4}{3}k,\frac{4}{3}(k+1))$ and consider the following statement: for all $s\in I_k$, there exists a non-decreasing function $C_s^{(\epsilon)}$ such that $\lVert u^{(\epsilon)}(t)\rVert_{H^s}\leq C_s(\| u_0 \|_{H^s})\langle t \rangle^{\alpha_k}$ for all $t \geq 0$ where $\alpha_k = \frac{1}{2}(3^{k+1}-1)$ for $k\geq 0$; note that $\alpha_k = 3\alpha_{k-1}+1$.

Let $k=0$ and fix $b\in (\frac{1}{2},1-\gamma), \gamma \in [\frac{1}{3},\frac{1}{2})$. By the local theory in $L^2(\mathbb{R})$ with $T=\Tilde{c} \lVert g\rVert_{L^2}^{-\rho}$ for some $\rho = \rho(b,\gamma)>0$ and $\Tilde{c}>0$,
\begin{equation*}
    \lVert u^{(\epsilon)}(t)\rVert_{H^s}\leq \lVert u_0 \rVert_{H^s} + C\lVert u^{(\epsilon)} \rVert_{X^{0,b}_{T}}^3\leq \lVert u_0 \rVert_{H^s}+C^\prime \lVert u_0 \rVert_{L^2}^3, t \in (0,T],
\end{equation*}
by \cref{mainprop}. Time evolving $u^{(\epsilon)}$ iteratively in $L^2(\mathbb{R})$ for $j=1,2,...$ and $t \in ((j-1)T,jT]$,
\begin{equation*}
\begin{split}
    \lVert u^{(\epsilon)}(t)\rVert_{H^s}&\leq \lVert u_0 \rVert_{H^s}+C^\prime j \lVert u_0 \rVert_{L^2}^3< \lVert u_0 \rVert_{H^s}+C^\prime (1+\frac{t}{T})\lVert u_0 \rVert_{L^2}^3\\
    &\leq
    \lVert u_0 \rVert_{H^s}+C^\prime (1+\frac{t}{T})\lVert u_0 \rVert_{H^s}^3
    \leq\lVert u_0 \rVert_{H^s}+C^\prime\lVert u_0 \rVert_{H^s}^3 +\frac{C^{\prime}}{\Tilde{c}}\lVert u_0 \rVert_{H^s}^{3+\rho}t.    
\end{split}
\end{equation*}
Noting that $\langle t \rangle \simeq 1+|t|$ and that $C^\prime, \Tilde{c}$ are constants that only depend on the given parameters, there exists $C_s^{(\epsilon)}$ such that
\begin{equation*}
    \lVert u_0 \rVert_{H^s}+C^\prime\lVert u_0 \rVert_{H^s}^3 +\frac{C^{\prime}}{\Tilde{c}}\lVert u_0 \rVert_{H^s}^{3+\rho}t \leq
    C_s^{(\epsilon)}(\lVert u_0 \rVert_{H^s})\langle t \rangle.
\end{equation*}
For example, $C_s^{(\epsilon)}(\zeta) =\frac{C^{\prime}}{\Tilde{c}}\zeta^{3+\rho}+ C^\prime\zeta^3+\zeta$ works.

For $k\geq 1$, assume the inductive hypothesis holds for $j=0,1,...,k-1$, and let $s \in I_k$. Fix $a\in (0,\frac{4}{3})$ such that $s-ja \in I_{k-j}$ for $j=1,2,...,k$ and $s-(k+1)a\leq 0$. Iteratively applying \cref{mainprop} and using the triangle inequaltiy $(a+b)^3\lesssim a^3+b^3$ for $a,b\geq 0$,
\begin{equation*}
\begin{split}
t \in (0,T]: \lVert u^{(\epsilon)}(t) \rVert_{H^s}&\leq \lVert u_0 \rVert_{H^s}+C(\lVert u_0 \rVert_{H^{s-a}}^3+\dots+\lVert u_0 \rVert_{H^{s-ka}}^{3^k})+C\lVert u_0 \rVert_{L^2}^{3^{k+1}}\\
    t\in(T,2T]:\lVert u^{(\epsilon)}(t)\rVert_{H^s}&\leq \lVert u_0 \rVert_{H^s}
    +C(\lVert u_0 \rVert_{H^{s-a}}^3+\dots+\lVert u_0 \rVert_{H^{s-ka}}^{3^k})+2C\lVert u_0 \rVert_{L^2}^{3^{k+1}}
    \\+C\Big(C_{s-a}^{(\epsilon)}&(\lVert u_0 \rVert_{H^{s-a}})^3\langle T\rangle^{3\alpha_{k-1}}+\dots+C_{s-ka}^{(\epsilon)}(\lVert u_0 \rVert_{H^{s-ka}})^{3^k}\langle T\rangle^{3^k\alpha_{0}}\Big)\\
    &\leq \lVert u_0 \rVert_{H^s}
    +C(\lVert u_0 \rVert_{H^{s-a}}^3+\dots+\lVert u_0 \rVert_{H^{s-ka}}^{3^k})+2C\lVert u_0 \rVert_{L^2}^{3^{k+1}}
    \\&+C\Big(C_{s-a}^{(\epsilon)}(\lVert u_0 \rVert_{H^{s}})^3+\dots+C_{s-ka}^{(\epsilon)}(\lVert u_0 \rVert_{H^{s}})^{3^k}\Big)\langle T\rangle^{3\alpha_{k-1}}\\
    t\in((j-1)T,jT]:\lVert u^{(\epsilon)}(t)\rVert_{H^s}&\leq \lVert u_0 \rVert_{H^s}
    +C(\lVert u_0 \rVert_{H^{s-a}}^3+\dots+\lVert u_0 \rVert_{H^{s-ka}}^{3^k})+jC\lVert u_0 \rVert_{L^2}^{3^{k+1}}
    \\
    +(j-1)C\Big(C_{s-a}^{(\epsilon)}&(\lVert u_0 \rVert_{H^{s}})^3+\dots+C_{s-ka}^{(\epsilon)}(\lVert u_0 \rVert_{H^{s}})^{3^k}\Big)\langle (j-1)T\rangle^{3\alpha_{k-1}}.    
\end{split}
\end{equation*}
Hence for all $t\geq 0$,
\begin{equation*}
\begin{split}
    \lVert u^{(\epsilon)}(t)\rVert_{H^s}&<\lVert u_0 \rVert_{H^s}
    +C(\lVert u_0 \rVert_{H^{s}}^3+\dots+\lVert u_0 \rVert_{H^{s}}^{3^k})+C(1+\frac{\lVert u_0 \rVert_{H^s}^\sigma}{\Tilde{c}}t)\lVert u_0 \rVert_{H^s}^{3^{k+1}}
    \\
    &+C\frac{\lVert u_0 \rVert_{H^s}^\sigma}{\Tilde{c}}\Big(C_{s-a}^{(\epsilon)}(\lVert u_0 \rVert_{H^{s}})^3+\dots+C_{s-ka}^{(\epsilon)}(\lVert u_0 \rVert_{H^{s}})^{3^k}\Big)\langle t\rangle^{3\alpha_{k-1}+1},    
\end{split}
\end{equation*}
which proves the claim.
\end{proof}
\begin{remark}
For $M=\mathbb{T}^d$, a simple example is constructed to show the failure of uniform well-posedness at $s<0$. Let $u_{n,k}(x,0)\coloneqq k \langle n\rangle^{-s}e^{in\cdot x}$ where $k>0,\: n \in \mathbb{Z}^d,\: x \in \mathbb{T}^d$ and $\cdot$ is the usual dot product. Let $\left\{k_n\right\}$ be a positive sequence that converges to $k$. By direct computation,
\begin{equation*}
    \lVert u_{n,k}(x,0)\rVert_{H^s}\simeq k,\: \lVert u_{n,k}(x,0)-u_{n,k_n}(x,0)\rVert_{H^s}\simeq |k_n-k| \xrightarrow[n\rightarrow\infty]{}0.
\end{equation*}
An exact solution to \eqref{adiabatic} corresponding to data $u_{n,k}(x,0)$ is given by
\begin{equation*}
    u_{n,k}(x,t) = k\langle n \rangle^{-s}e^{-it(\epsilon^2 |n|^4+|n|^2-k^2\langle n \rangle^{-2s})+in\cdot x}.
\end{equation*}
Given any $T>0$ and $t\in [0,T]$,
\begin{equation}\label{illest}
\begin{split}
    \lVert u_{n,k_n}(x,t)-u_{n,k}(x,t)\rVert_{H^s}&\simeq_s |k_ne^{itk_n^2\langle n\rangle^{-2s}}-ke^{itk^2\langle n\rangle^{-2s}}|\\
    &\geq |k|\cdot|e^{it\langle n \rangle^{-2s}(k_n^2-k^2)}-1|-|k_n-k|.    
\end{split}
\end{equation}
By defining $k_n\coloneqq \Big(k^2+t^{-1}\langle n \rangle^{2s}\pi\Big)^{1/2}$, which is shown to converge to $k$ since $s<0$, the RHS of \eqref{illest} is bounded below by $k$ for all but finitely many $n$.
\end{remark}

Smoothing due to dispersion depends on the spatial domain, and therefore our analysis on $\mathbb{T}$ needs a modification. We consider the continuous embedding of $X^{s,b}$ into $L^p(\mathbb{T}\times\mathbb{R})$ for $p=4,6$ where Fourier-analytic tools are available. This method is motivated from \cite{oh2017quasi} where it is shown that
\begin{equation*}
    \lVert u \rVert_{L^4(\mathbb{T}\times\mathbb{R})}\lesssim \lVert u \rVert_{X^{0,\frac{5}{16}}},
\end{equation*}
with $w(k) = k^4$.

\begin{proposition}\label{thm4.1}
For $w(k)=\epsilon^2k^4+k^2$ with $\epsilon>0$,
\begin{equation}\label{strichartz}
\lVert u \rVert_{L^4(\mathbb{T}\times\mathbb{R})}\lesssim \epsilon^{-\frac{1}{8}}\lVert u \rVert_{X^{0,\frac{5}{16}}},\:\lVert u \rVert_{L^6(\mathbb{T}\times\mathbb{R})}\lesssim \epsilon^{-\frac{1}{6}}\lVert u \rVert_{X^{0,\frac{5}{12}}}.
\end{equation}
\end{proposition}
The end of this section contains an example that shows the sharpness of \cref{thm4.1}.
\begin{remark}\label{rmktobeproved}
Estimates of the form \eqref{strichartz} can be generalized to higher order derivatives. For instance if $w(k)=k^\delta$ with $\delta\in\left\{2,3,4,\dots\right\}$, then $\lVert u \rVert_{L^4(\mathbb{T}\times\mathbb{R})}\lesssim \lVert u \rVert_{X^{0,\frac{\delta+1}{4\delta}}}$. See the end of this section. 
\end{remark}

When $\epsilon=0$, however, the embedding $X^{0,b}\hookrightarrow L^6(\mathbb{T}\times\mathbb{R})$ for $b<\frac{1}{2}$ fails to hold. If it were to hold, then NLSE with the nonlinearity $\pm|u|^4u$ is well-posed in $L^2(\mathbb{T})$ by the contraction mapping argument, which contradicts \cite[Corollary 1.3]{kishimoto2014remark} where Kishimoto showed that the data-to-solution map, if it exists, cannot be $C^5$. For the quintic nonlinearity, our method yields a sequence of solutions, say $\left\{u^{(\epsilon)}\right\}_{\epsilon>0}$, with $u^{(\epsilon)}(0)=u_0\in L^2(\mathbb{T})$. For $T>0,\:p\in (1,\infty)$, since $\left\{u^{(\epsilon)}\right\}_{\epsilon>0}\subseteq C([0,T],L^2(\mathbb{T}))\hookrightarrow L^p([0,T],L^2(\mathbb{T}))$ with $\lVert u^{(\epsilon)} \rVert_{L^p([0,T],L^2(\mathbb{T}))} = T^{\frac{1}{p}}\lVert u_0 \rVert_{L^2}$, there exists a weakly convergent subsequence in $L^p([0,T],L^2(\mathbb{T}))$. However without any extra regularity on the solutions, it is insufficient to show that the limit defines a strong solution. Hence our method of adding a small fourth-order dispersion seems unlikely to yield positive results for the mass-critical periodic NLSE.

\begin{proof}[Proof of \cref{thm311} for $M=\mathbb{T}$]
For $s\geq 0$, $\frac{1}{2}<b<\frac{11}{16}$, and $T\leq 1$, let $X,\Gamma$ be as in \eqref{functionspace}, \eqref{contraction}, respectively. To show that $\Gamma$ is a contraction, we use the following lemma.
\begin{lemma}\label{lemmaembedding}
For every $s\geq 0$ and $\epsilon>0$,
\begin{equation*}
\begin{split}
\lVert J_{\epsilon}(f\overline{g})h\rVert_{X^{s,-\frac{5}{16}}}\lesssim _{s,\epsilon} \lVert  f \rVert_{X^{s,\frac{5}{16}}}\lVert g \rVert_{X^{0,\frac{5}{16}}}\lVert h \rVert_{X^{0,\frac{5}{16}}}+\lVert  f \rVert_{X^{0,\frac{5}{16}}}\lVert g \rVert_{X^{s,\frac{5}{16}}}\lVert h \rVert_{X^{0,\frac{5}{16}}}\\
    +\lVert  f \rVert_{X^{0,\frac{5}{16}}}\lVert g \rVert_{X^{0,\frac{5}{16}}}\lVert h \rVert_{X^{s,\frac{5}{16}}}.
\end{split}    
\end{equation*}
\end{lemma}
\begin{proof}[Proof of lemma \ref{lemmaembedding}]
By duality if $\overline{v}\in (X^{s,-\frac{5}{16}})^*$,
\begin{align*}
    |\langle J_\epsilon(f \overline{g})h,v \rangle|_{L^2_{x,t}} &= |\langle \langle \nabla \rangle^s (J_\epsilon(f \overline{g})h),\langle \nabla \rangle^{-s}v\rangle_{L^2}|\leq \lVert \langle \nabla \rangle^s (J_\epsilon(f\overline{g})h) \rVert_{L^{4/3}} \lVert \langle \nabla\rangle^{-s}v \rVert_{L^4}\\
    &\lesssim
\lVert \langle \nabla \rangle^s (J_\epsilon(f\overline{g})h) \rVert_{L^{4/3}} \lVert v \rVert_{X^{-s,\frac{5}{16}}},
\end{align*}
where the last inequality is by \cref{thm4.1}. The desired estimate follows by repeatedly applying the Leibniz rule for Sobolev spaces on $\lVert \langle \nabla \rangle^s (J_\epsilon(f\overline{g})h) \rVert_{L^{4/3}}$.
\end{proof}
By \cref{lemmaembedding,l:32},
\begin{equation}\label{globaliterative}
\begin{split}
\lVert \Gamma u\rVert_{X^{s,b}_T}
    &\lesssim \lVert u_0 \rVert_{H^s}+\lVert J_\epsilon(|u|^2)u\rVert_{X^{s,b-1}_T}\lesssim \lVert u_0\rVert_{H^s}+T^{\frac{11}{16}-b} \lVert J_\epsilon(|u|^2)u \rVert_{X^{s,-\frac{5}{16}}_T}\\
    &\lesssim \lVert u_0\rVert_{H^s}+ T^{(\frac{5}{4}-b)-}\lVert u \rVert_{X^{0,b}_T}^2 \lVert u\rVert_{X^{s,b}_T}\\
    &\lesssim \lVert u_0\rVert_{H^s}+T^{(\frac{5}{4}-b)-}\lVert u \rVert_{X^{s,b}_T}^3.
\end{split}    
\end{equation}
Hence $\Gamma$ is well-defined on some closed ball of $X^{s,b}_T$ of radius $C\lVert u_0\rVert_{H^s}$ where $C$ is sufficiently large. The difference $\lVert \Gamma u-\Gamma v\rVert_{X^{s,b}_T}$ is estimated similarly. Uniqueness extends to the full $X^{s,b}_T$ by the standard continuity argument. The local Lipschitz regularity of the data-to-solution map follows from a chain estimates similar to \eqref{globaliterative}. 

To extend the solution $u\in C([0,T],H^s(\mathbb{T}))\cap X^{s,b}_T$ globally in time, the global well-posedness at $L^2$ (by the mass conservation) is used to the estimate \eqref{globaliterative} to obtain $T_0 = T_0(\lVert u_0\rVert_{L^2})>0$ such that $\lVert u \rVert_{X^{s,b}_{T_0}}\lesssim \lVert u_0\rVert_{H^s}$. Since $b >\frac{1}{2}$, there exists some constant $C>0$ such that $\lVert u(t)\rVert_{H^s}\leq C \lVert u_0\rVert_{H^s}$ on $t\in [0,T_0]$ by \cref{l:32}. Iterating this procedure,
\begin{equation*}
\lVert u(t)\rVert_{H^s}\leq C^{|t|}\lVert u_0\rVert_{H^s}  
\end{equation*}
holds for all $t \in \mathbb{R}$.
\end{proof}
\begin{proof}[Proof of \cref{thm4.1}]
We closely follow the method in \cite[Proposition 2.13]{tao2006nonlinear}. Let $\| u \|_{L^p}:= \| u \|_{l^p_k L^p_\tau}$. For $m\in \mathbb{N}\cup\left\{0\right\}$, define the dyadic projector $\chi_m (\tau,k) \coloneqq \chi_{2^m\leq \langle \tau + w(k)\rangle<2^{m+1}}$ and $\widehat{u_{2^m}}(\tau,k) \coloneqq \widehat{u}(\tau,k)\chi_m(\tau,k)$. Then
\begin{equation*}
    \lVert u \rVert_{L^4}^2 = \lVert u^2 \rVert_{L^2}\lesssim \sum_{m,n\geq 0} \lVert u_{2^m}u_{2^{m+n}}\rVert_{L^2} = \sum_{m,n\geq 0}\lVert \widehat{u_{2^m}}\ast \widehat{u_{2^{m+n}}}\rVert_{L^2}.
\end{equation*}
Since $\chi_m^2 = \chi_m$,
\begin{equation*}
\begin{split}
     &\widehat{u_{2^m}}\ast \widehat{u_{2^{m+n}}}(\tau,k)\\
     = &\sum_{k_1\in\mathbb{Z}}\int \widehat{u_{2^m}}(\tau_1,k_1)\widehat{u_{2^{m+n}}}(\tau-\tau_1,k-k_1)\chi_m(\tau_1,k_1)\chi_{m+n}(\tau-\tau_1,k-k_1)d\tau_1.    
\end{split}
\end{equation*}
By the Cauchy-Schwarz inequality in $\tau_1,k_1$, H\"older's inequality in $l_k^2 L_\tau^\infty$, and Young's inequality,
\begin{equation*}
\begin{split}
    \lVert \widehat{u_{2^m}}\ast \widehat{u_{2^{m+n}}}\rVert_{L^2}\leq \lVert \chi_m\ast \chi_{m+n}\rVert_{l^\infty_{k}L^\infty_\tau}^{1/2} \left|\left| |\widehat{u_{2^m}}|^2\ast |\widehat{u_{2^{m+n}}}|^2\right|\right|_{l^1_{k}L^1_\tau}^{1/2}\\
    \leq \lVert \chi_m\ast \chi_{m+n}\rVert_{l^\infty_{k}L^\infty_\tau}^{1/2} \lVert u_{2^m} \rVert_{L^2}\lVert u_{2^{m+n}} \rVert_{L^2}.    
\end{split}
\end{equation*}
It remains to extract a sufficient decay (in $m,n$) from $\lVert \chi_m\ast \chi_{m+n}\rVert_{L^\infty}^{1/2}$. Since
\begin{equation*}
    \chi_m\ast \chi_{m+n}(\tau,k) = \sum_{k_1} \int \chi_{m}(\tau_1,k_1)\chi_{m+n}(\tau-\tau_1,k-k_1)d\tau_1,
\end{equation*}
each non-zero integral, for a fixed $k_1$, is at most $O(2^m)$. It remains to count how many $k_1$s gives rise to non-zero integrals. For a fixed $\tau,k$, if the integral is non-zero, then
\begin{equation*}
    \tau_1 + w(k_1) = O(2^m),\: \tau-\tau_1 + w(k-k_1) = O(2^{m+n}),
\end{equation*}
and therefore
\begin{equation}\label{fouriercomb}
    \tau + w(k_1)+w(k-k_1) = O(2^{m+n}). 
\end{equation}
It can be shown by direct computation that $x\mapsto w(x) + w(k-x) = w(x) + w(x-k)$ has a global minimum at $x=\frac{k}{2}$. Changing variable $k_1^\prime = k_1 - \frac{k}{2}$ so that in the new coordinate, $x\mapsto w(x)+w(k-x)$ is centered at the origin, the LHS of \eqref{fouriercomb} becomes
\begin{equation*}
\begin{split}
    &\tau + 2\epsilon^2 k_1^4 + (2+3\epsilon^2 k^2)k_1^2 + \frac{k^2(4+\epsilon^2k^2)}{8}\\ = &\tau + 2\epsilon^2 \Big(k_1^2 + \frac{1}{2}(\epsilon^{-2}+\frac{3k^2}{2})\Big)^2 - \Big(\epsilon^2k^4 + k^2 + \frac{1}{2\epsilon^2}\Big),    
\end{split}
\end{equation*}
where $k_1^\prime$ is re-labelled to $k_1$. This implies that $k_1$s are constrained in finitely many intervals of length $O(\epsilon^{-\frac{1}{2}}2^{\frac{m+n}{4}})$, and therefore
\begin{equation*}
    \lVert \chi_m\ast \chi_{m+n}\rVert_{L^\infty}^{1/2} \lesssim \epsilon^{-\frac{1}{4}}2^{\frac{5m}{8}+\frac{n}{8}},
\end{equation*}
from which
\begin{equation}\label{fouriercomb2}
\begin{split}
    \lVert u \rVert_{L^4}^2\lesssim \epsilon^{-\frac{1}{4}}\sum_{m,n\geq 0} 2^{\frac{5m}{8}+\frac{n}{8}}\lVert u_{2^m}\rVert_{L^2}\lVert u_{2^{m+n}}\rVert_{L^2}\\
    = \epsilon^{-\frac{1}{4}}\sum_{n\geq 0} 2^{-\frac{3n}{16}}\sum_{m\geq 0} 2^{\frac{5m}{16}}\lVert u_{2^m}\rVert_{L^2}2^{\frac{5}{16}(m+n)}\lVert u_{2^{m+n}} \rVert_{L^2}\lesssim \epsilon^{-\frac{1}{4}} \lVert u \rVert_{X^{0,\frac{5}{16}}}^2,    
\end{split}
\end{equation}
where the Cauchy-Schwarz inequality is used at the last inequality and 
\begin{equation*}
\lVert u \rVert_{X^{0,b}}^2 \simeq \sum\limits_{m\geq 0} 2^{2bm}\lVert u_{2^m}\rVert_{L^2}^2,    
\end{equation*}
by the Plancherel's Theorem.

Extending the previous counting method, the $L^6$ estimate is obtained. Observe that
\begin{equation*}
\begin{split}
\lVert u \rVert_{L^6}^3 &= \lVert u^3\rVert_{L^2}\lesssim \sum_{m,n,l \geq 0} \lVert u_{2^m}u_{2^{m+n}}u_{2^{m+n+l}}\rVert_{L^2} = \sum_{m,n,l\geq 0} \lVert \widehat{u_{2^{m}}}\ast\widehat{u_{2^{m+n}}}\ast\widehat{u_{2^{m+n+l}}}\rVert_{L^2}\\
    =\sum_{m,n,l\geq 0}& \left|\left|\sum_{k_1,k_2\in\mathbb{Z}}\iint_{\tau_1,\tau_2} \widehat{u_{2^{m}}}(\tau_1,k_1)\widehat{u_{2^{m+n}}}(\tau_2,k_2)\widehat{u_{2^{m+n+l}}}(\tau-\tau_1-\tau_2,k-k_1-k_2) \right|\right|_{L^2}\\
    &\leq \sum_{m,n,l\geq 0} \lVert \chi_m\ast\chi_{m+n}\ast\chi_{m+n+l}\rVert_{L^\infty_{k,\tau}}^{1/2}\lVert u_{2^m} \rVert_{L^2_{x,t}}\lVert u_{2^{m+n}} \rVert_{L^2_{x,t}}\lVert u_{2^{m+n+l}} \rVert_{L^2_{x,t}}.   
\end{split}    
\end{equation*}
We claim $\lVert \chi_m\ast\chi_{m+n}\ast\chi_{m+n+l}\rVert_{L^\infty_{k,\tau}}^{1/2}\lesssim \epsilon^{-\frac{1}{2}}2^{\frac{5m+3n+l}{4}}$. Using the support conditions of dyadic projectors,
\begin{equation*}
\begin{split}
\tau_1+w(k_1) = O(2^m),\: \tau_2 + w(k_2) &= O(2^{m+n})\\
    \tau-\tau_1-\tau_2 + w(k-k_1-k_2) &= O(2^{m+n+l}),  
\end{split}
\end{equation*}
and therefore
\begin{align*}
    \tau + w(k_1)+w(k_2)+w(k-k_1-k_2) &= O(2^{m+n+l}).
\end{align*}
Define 
\begin{align*}
    E(\tau,k) = \left\{(k_1,k_2)\in\mathbb{Z}^2:\tau + w(k_1)+w(k_2)+w(k-k_1-k_2) = O(2^{m+n+l})\right\}.
\end{align*}
The map $(k_1,k_2) \mapsto \Tilde{w}(k_1,k_2)\coloneqq w(k_1)+w(k_2)+w(k-k_1-k_2)$, for every fixed $k\in\mathbb{Z}$, has a global minimum at $(\frac{k}{3},\frac{k}{3})$ by direct computation. We change variables $(k_1,k_2)\mapsto (k_1+\frac{k}{3},k_2+ \frac{k}{3})$, after which $\Tilde{w}$ in the new variables, without re-labelling, is a polynomial of degree 4 in two variables centered at the origin. Transform $\Tilde{w}$ further by considering another change of variable, $k_1 = r\cos(\theta),\: k_2 = r\sin(\theta)$, and considering $\Tilde{w}$ as a polynomial of one variable, namely $r$, for every fixed $\theta \in [0,2\pi)$. Then,
\begin{equation}
\begin{split}
v(r)\coloneqq\Tilde{w}(r,\theta) &= \frac{r^2}{12}\bigg( 3\epsilon^2 \Big(9-\cos(4\theta)+8\sin(2\theta)\Big)r^2\\
    &-12\epsilon^2 k\Big(\cos(\theta)-\cos(3\theta)+\sin(\theta)+\sin(3\theta)\Big)r\\
    &\hspace{3cm}+(8\epsilon^2k^2+12)(\sin(2\theta)+2)\bigg)
    +\frac{k^2}{27}(\epsilon^2k^2+9).    
\end{split}    
\end{equation}
We claim $v$ is an increasing and convex function on $r>0$ for every $\theta\in[0,2\pi),\: k \in \mathbb{Z}$. Since $v$ has no $r^1$-term, it is clear that $v^\prime(0)=0$, and it is shown by direct computation that $v^{\prime\prime}(r)\geq 0$ on $r\geq 0$. Indeed
\begin{equation*}
\begin{split}
v^{\prime\prime}(r)=6\epsilon^2 \Big(\sin(2\theta)&+2\Big)^2r^2 +6\epsilon^2 k \Big( \cos (3 \theta)- (\sin (\theta)+\sin (3 \theta)+\cos (\theta))\Big)r\\
&+\Big(\frac{4}{3} \epsilon ^2 k^2+2\Big)\Big(\sin (2 \theta)+2\Big),
\end{split} 
\end{equation*}
which, by another direct computation, is non-negative for every $\theta\in [0,2\pi),\:k\in\mathbb{Z}$. Hence for every $k\in\mathbb{Z}$, $|E(\tau,k)|$ is maximized at $\tau = -\frac{k^2}{27}(\epsilon^2k^2+9)$ since the $\tau$-term corresponds to translating $\Tilde{w}$. To obtain a uniform estimate in $k$, note that the following lower bound
\begin{equation*}
    \frac{1}{r^2}\Big(v(r)-\frac{k^2}{27}(\epsilon^2k^2+9)\Big)\gtrsim\epsilon^2 r^2
\end{equation*}
holds for all $r\geq 0$ with the implicit constant independent of $\theta,k$. Then
\begin{equation*}
    \sup_{\tau\in\mathbb{R}}|E(\tau,k)| \lesssim \left|\left\{(k_1,k_2)\in\mathbb{Z}^2: \epsilon^2 (k_1^2+k_2^2)^2 = O(2^{m+n+l})\right\}\right|\lesssim \epsilon^{-1}2^{\frac{m+n+l}{2}},
\end{equation*}
and together with the $O(2^{2m+n})$ contribution from each integral corresponding to $(k_1,k_2)\in E(\tau,k)$, we have
\begin{equation*}
    \lVert \chi_{m}\ast\chi_{m+n}\ast \chi_{m+n+l}\rVert_{L^\infty_{k,\tau}}\lesssim \epsilon^{-1}2^{\frac{m+n+l}{2}+2m+n} = \epsilon^{-1}2^{\frac{5m+3n+l}{2}}.
\end{equation*}
The rest follows immediately as \eqref{fouriercomb2}.
\end{proof}
\begin{proof}[Proof of remark \ref{rmktobeproved}]
We estimate $\lVert \chi_m \ast \chi_{m+n}\rVert_{L^\infty_{k,\tau}}^{1/2}$ in $m$ and $n$ such that the sum
\begin{equation*}
\sum\limits_{m,n\geq 0}\lVert \chi_m \ast \chi_{m+n}\rVert_{L^\infty_{k,\tau}}^{1/2} \lVert u_{2^m}\rVert_{L^2}\lVert u_{2^{m+n}}\rVert_{L^2}    
\end{equation*} 
converges. As in the proof of \cref{thm4.1}, the problem is reduced to estimating the cardinality of
\begin{equation*}
    E(\tau,k) = \left\{k_1\in\mathbb{Z}:\tau + w(k_1)+w(k-k_1) = O(2^{m+n})\right\}.
\end{equation*}
Assume $w(k)=k^\delta$ for $\delta\geq 2$ even. Observe that $x\mapsto w(x)+w(k-x)$ is convex with a global minimum at $x=\frac{k}{2}$. Changing variable $x\mapsto x+\frac{k}{2}$, we have $|E(\tau,k)| = |\Tilde{E}(\tau,k)|$ where
\begin{equation*}
    \Tilde{E}(\tau,k) = \left\{k_1 \in \mathbb{Z}: \tau + w(k_1+\frac{k}{2})+w(\frac{k}{2}-k_1) = O(2^{m+n})\right\}.
\end{equation*}
For a fixed $k\in \mathbb{Z}$, since $\Tilde{w}(k_1) \coloneqq w(k_1+\frac{k}{2})+w(\frac{k}{2}-k_1)$ is even (in $k_1$) and convex with a global minimum of $2\cdot (\frac{k}{2})^\delta$, and the $\tau$-term corresponds to translating $\Tilde{w}(k_1)$,
\begin{equation*}
\begin{split}
\sup_{\tau \in \mathbb{R}}|\Tilde{E}(\tau,k)|&\lesssim \left|\left\{k_1: \Tilde{w}(k_1)-2\cdot \Big(\frac{k}{2}\Big)^\delta = O(2^{m+n})\right\} \right|\\
    &\leq \left|\left\{k_1: k_1^\delta = O(2^{m+n})\right\}\right|\lesssim 2^{\frac{m+n}{\delta}},
\end{split}
\end{equation*}
where the second inequality holds since $\Tilde{w}(k_1)-2\cdot \Big(\frac{k}{2}\Big)^\delta\geq k_1^\delta$ for all $k_1,k \in \mathbb{Z}$. Hence
\begin{equation*}
\lVert \chi_m \ast \chi_{m+n}\rVert_{L^\infty_{k,\tau}}^{1/2}\lesssim 2^{\frac{(\delta+1)m}{2\delta}+\frac{n}{2\delta}}    
\end{equation*}
and the argument proceeds as \eqref{fouriercomb2}.

Assume $\delta$ is odd. By the argument in \cite[Theorem 3.18]{erdougan2016dispersive}, 
\begin{equation*}
    \lVert \widehat{u_{2^m}}\ast \widehat{u_{2^{m+n}}} \rVert_{L^2_\tau l^2_{|k|\leq 2^a}}\lesssim 2^{\frac{a+m}{2}}\lVert u_{2^m}\rVert_{L^2}\lVert u_{2^{m+n}}\rVert_{L^2}
\end{equation*}
for $a \geq 0$ when $k = O(1)$.

To estimate $\lVert \chi_m\ast \chi_{m+n}\rVert_{L^\infty_\tau l^\infty_{|k|>2^a}}^{1/2}$, note that each non-zero integral has $O(2^m)$ contribution. After changing variable $k_1\mapsto k_1+\frac{k}{2}$, $w(k_1)+w(k-k_1)$ becomes $w(k_1+\frac{k}{2})-w(k_1-\frac{k}{2})$, which by direct computation amounts to
\begin{equation*}
\begin{split}
\Tilde{w}(k_1)&=(k_1+\frac{k}{2})^\delta - (k_1-\frac{k}{2})^\delta = 2\sum_{j\:\text{even},\:0\leq j \leq \delta} \binom{\delta}{j}k_1^j (\frac{k}{2})^{\delta-j}\\
    &= 2k\sum_{j\:\text{even},\:0\leq j \leq \delta} \binom{\delta}{j}k_1^j \frac{k^{\delta-j-1}}{2^{\delta-j}}.    
\end{split}
\end{equation*}
Since $\Tilde{w}$ is even and convex with the global minimum $2\cdot\Big(\frac{k}{2}\Big)^\delta$,
\begin{equation*}
\begin{split}
\sup\limits_{\tau \in \mathbb{R}}|\Tilde{E}(\tau,k)|&\lesssim \left|\left\{k_1: \Tilde{w}(k_1)-2\cdot \Big(\frac{k}{2}\Big)^\delta = O(2^{m+n})\right\} \right| \\
  &= \left|\left\{k_1:2k\sum_{j\:\text{even},\:2\leq j \leq \delta} \binom{\delta}{j}k_1^j \frac{k^{\delta-j-1}}{2^{\delta-j}} = O(2^{m+n})\right\}\right|\\
  &\lesssim \left|\left\{k_1:\sum_{j\:\text{even},\:2\leq j \leq \delta} \binom{\delta}{j}k_1^j \frac{k^{\delta-j-1}}{2^{\delta-j}} = O(2^{m+n-a})\right\}\right|\\
  &\leq \left|\left\{k_1: \frac{\delta}{2}k_1^{\delta-1}=O(2^{m+n-a})\right\}\right|,    
\end{split}  
\end{equation*}
which implies $\lVert \chi_m\ast \chi_{m+n}\rVert_{L^\infty_\tau l^\infty_{|k|>2^a}}^{1/2}\lesssim 2^{\frac{\delta m + n-a}{2(\delta-1)}}$. Set $a=\frac{m+n}{\delta}$ to equate the bounds from the low and high frequencies to obtain
\begin{equation*}
    \lVert u_{2^m}u_{2^{m+n}}\rVert_{L^2_{x,t}}\lesssim 2^{\frac{(\delta+1)m}{2\delta}+\frac{n}{2\delta}}\lVert u_{2^m}\rVert_{L^2}\lVert u_{2^{m+n}}\rVert_{L^2},
\end{equation*}
from which the argument proceeds as in \eqref{fouriercomb2}. 
\end{proof}
\begin{remark}
As in \cite[Footnote 9]{oh2017quasi}, we use projectors in the space-time Fourier space to show that \cref{thm4.1} is sharp. Without loss of generality, let $w(k)=k^\delta$ for $\delta \in \left\{2,3,\dots\right\}$. Define $\widehat{u_N}(k,\tau)\coloneqq \chi_{N}(k)\chi_{N^\delta}(\tau)$ where $\chi_{N}(k)$ is the characteristic function on $k \in [-N,N]$ and similarly for $\chi_{N^\delta}(\tau)$. A direct computation reveals
\begin{equation*}
    \lVert u_N \rVert_{L^4_{x,t}}\simeq N^{\frac{3}{4}(1+\delta)},\: \lVert u_N \rVert_{L^6_{x,t}} \simeq N^{\frac{5}{6}(1+\delta)},\: \lVert u_N \rVert_{X^{0,b}}\simeq N^{\frac{1+(2b+1)\delta}{2}},
\end{equation*}
and this shows sharpness. In fact by considering $p=2q$ for $q\in \left\{2,3,\dots\right\}$, we derive $\lVert u_N\rVert_{L^{2q}_{x,t}}\simeq N^{\frac{2q-1}{2q}(1+\delta)}$. Hence for $X^{0,b}\hookrightarrow L^{2q}(\mathbb{T}\times \mathbb{R})$ to hold, it is necessary that 
\begin{equation}\label{highorder}
b\geq \frac{(q-1)(1+\delta)}{2q\delta}. 
\end{equation}
From the proof of \cref{thm311}, the local well-posedness on $\mathbb{T}$ follows from the embedding of the form $X^{0,b}\hookrightarrow L^p(\mathbb{T}\times\mathbb{R})$ with $b<\frac{1}{2}$. Combining with \eqref{highorder}, we obtain $q-1<\delta$, and hence $q_{max}=\delta$, or equivalently $p_{max}=2\delta$. For $\delta \in \left\{2,3,\dots\right\}$, considering that the highest odd-power $L^2$-subcritical nonlinearity of the model
\begin{equation}\label{high order dispersive}
    \begin{cases}
    i\partial_t u + \partial_x^\delta u =\pm|u|^{\rho-1}u,\: (x,t)\in \mathbb{T}\times \mathbb{R}\\
    u(0)=u_0\in L^2(\mathbb{T}),
    \end{cases}
\end{equation}
is $\rho_{\text{sub}} = 2\delta-1$, it is observed that \eqref{high order dispersive} is globally well-posedness in $L^2(\mathbb{T})$ for all $\rho =2j-1,\: j=2,3,\dots,\delta$, provided $X^{0,b}\hookrightarrow L^{p_{max}}(\mathbb{T}\times\mathbb{R})$ holds.
\end{remark}
\section{Semi-classical Limit}\label{section limit}
\begin{proof}[Proof of \cref{thmII}]
We prove the claim only for $M=\mathbb{R}$ since $H^s(M)$ is an algebra for $s>\frac{1}{2}$, and therefore, substituting certain integrals by sums admits a similar proof for $M=\mathbb{T}$. 

Let $U_\epsilon(t) = e^{-it(\epsilon^2 \xi^4 + \xi^2)}$ and $U(t):= U_0(t)$. Let $\delta>0$. We claim there exists $\epsilon_0>0$ such that if $\epsilon\in (0,\epsilon_0)$, then  
\begin{equation*}
    \sup\limits_{t\in [0,T]}\lVert U_\epsilon(t)u_0^{(\epsilon)}-U(t)u_0\rVert_{H^1}<\delta.
\end{equation*}
Observe
\begin{equation*}
    \lVert U_\epsilon(t)u_0^{(\epsilon)}-U(t)u_0\rVert_{H^1}\leq \lVert u_0^{(\epsilon)}-u_0\rVert_{H^1}+\lVert(U_\epsilon(t)-U(t))u_0\rVert_{H^1}.
\end{equation*}
Since $|e^{-i\epsilon^2t |\xi|^4}-1|^2=2(1-\cos(\epsilon^2t|\xi|^4))$, it is shown that $\lim\limits_{\epsilon \rightarrow 0+}I=0$ where $I=\int (1-\cos(\epsilon^2t|\xi|^4)) \langle \xi \rangle^{2}|\widehat{u_0}|^2d\xi$. Let $\xi_n = \xi_n(\epsilon,t)$ be the first positive root to $1-\cos(\epsilon^2t\xi^4) = \frac{1}{n}$. Estimating $I$ in two different regions,
\begin{equation*}
\begin{split}
    I &= \int_{|\xi|\leq \xi_n(\epsilon,t)} + \int_{|\xi|>\xi_n(\epsilon,t)} \leq \frac{\lVert u_0 \rVert_{H^1}^2}{n}+2\int_{|\xi|>\xi_n(\epsilon,T)} \langle\xi\rangle^{2}|\widehat{u_0}|^2d\xi,   
\end{split}
\end{equation*}
which yields the desired result by taking $n\rightarrow\infty,\epsilon\rightarrow 0$.

To proceed with the nonlinear estimates, observe
\begin{equation*}
\begin{split}
\int_0^t U_\epsilon(t-t^\prime)\Big(J_\epsilon(|u^{(\epsilon)}|^2)u^{(\epsilon)}\Big)(t^\prime)-U(t-t^\prime)\Big(|u|^2u\Big)(t^\prime)dt^\prime = I_1&+I_2+I_3+I_4,\:\text{where}\\
I_1 = \int_0^t U_\epsilon(t-t^\prime)\Big( J_\epsilon(|u^{(\epsilon)}|^2)(u^{(\epsilon)}-u)\Big)dt^\prime;\:I_2=\int_0^t U_\epsilon(t-t^\prime)&\Big(J_\epsilon\Big(|u^{(\epsilon)}|^2-|u|^2\Big)u\Big)dt^\prime.\\
I_3 = \int_0^t U_\epsilon(t-t^\prime)\Big(\big(J_\epsilon(|u|^2)-|u|^2\big)u\Big)dt^\prime;\: I_4= \int_0^t \Big(U_\epsilon(t-t^\prime&)-U(t-t^\prime)\Big)(|u|^2u)dt^\prime.    
\end{split}
\end{equation*}
By the triangle inequality and \cref{cazenave},
\begin{equation*}
\begin{split}
    \lVert I_1\rVert_{H^1}&\leq \int_0^t \lVert J_\epsilon(|u^{(\epsilon)}|^2)(u^{(\epsilon)}-u)\rVert_{H^1}dt^\prime \lesssim \Big(\sup\limits_{\epsilon>0}\lVert u^{(\epsilon)}\rVert_{C_T H^1}\Big)^2\int_0^t \lVert u^{(\epsilon)}-u\rVert_{H^1}dt^\prime\\
    \lVert I_2\rVert_{H^1}&\lesssim (\sup\limits_{\epsilon>0}\lVert u^{(\epsilon)}\rVert_{C_TH^1}+\lVert u \rVert_{C_T H^1})\lVert u \rVert_{C_TH^1}\int_0^t \lVert u^{(\epsilon)}-u\rVert_{H^1}dt^\prime,    
\end{split}
\end{equation*}
where $\| u \|_{C_T H^1} \leq C(\| u_0 \|_{H^1})$ by the energy conservation of NLSE and $\sup\limits_{\epsilon>0}\lVert u^{(\epsilon)}\rVert_{C_TH^1} \leq C$ by the boundedness of $u_0^{(\epsilon)}$ in $H^2$. To estimate $I_3$, note that
\begin{equation*}
  \lVert I_3\rVert_{H^1}\lesssim \lVert u \rVert_{C^0_TH^1}\int_0^t \lVert (J_\epsilon-I)(|u|^2)\rVert_{H^1}dt^\prime.  
\end{equation*}
By the Dominated Convergence Theorem,
\begin{equation*}
F_\epsilon(t):=\int_0^t \lVert (J_\epsilon-I)(|u|^2)\rVert_{H^1}dt^\prime\xrightarrow[\epsilon\rightarrow 0]{}0,    
\end{equation*}
for every $t \in [0,T]$. Since $|F_\epsilon(t)-F_\epsilon(t^\prime)|\lesssim\lVert u \rVert_{C_TH^1}^2 |t-t^\prime|$, the family $\left\{F_\epsilon\right\}\subseteq C([0,T];\mathbb{R})$ is uniformly equicontinuous, and thus $\| I_3 \|_{H^1}\xrightarrow[\epsilon\rightarrow 0]{}0$ uniformly in $t$ by the Arzel\`a-Ascoli Theorem.

To estimate $I_4$, a similar Arzel\`a-Ascoli argument is used. More precisely, we have
\begin{equation*}
    \lVert I_4 \rVert_{H^1}\leq \int_0^t \left|\left|\Big(U_\epsilon(t-t^\prime)-U(t-t^\prime)\Big)(|u|^2u)\right|\right|_{H^1} dt^\prime = \int_0^t \left|\left|\Big(U_\epsilon(t^\prime)-U(t^\prime)\Big)(|u(t-t^\prime)|^2u(t-t^\prime))\right|\right|_{H^1} dt^\prime,
\end{equation*}
and the last term tends to zero uniformly in $t$ as $\epsilon\rightarrow 0$.

Hence for all $\epsilon>0$ sufficiently small,
\begin{equation*}
    \lVert u^{(\epsilon)}(t)-u(t)\rVert_{H^1}\lesssim \delta\langle t\rangle +\int_0^t \lVert u^{(\epsilon)}(t^\prime)-u(t^\prime)\rVert_{H^1}dt^\prime,
\end{equation*}
from which the Gronwall's inequality yields
\begin{equation*}
    \lVert u^{(\epsilon)}(t)-u(t)\rVert_{H^1}\leq C\delta \langle t \rangle e^{Ct},
\end{equation*}
for all $t \in [0,T]$ for some $C>0$. Since $\delta>0$ is arbitrary, the proof is complete.
\end{proof}
\begin{remark}
Under similar hypotheses as \cref{thmII}, the continuity of solutions to \eqref{adiabatic} in $\epsilon \in [0,\infty)$ can be shown on a compact time interval by a similar Gronwall's inequality argument.
\end{remark}
\begin{remark}
The proof of \cref{thmII} cannot be extended to $T=\infty$. Examples are provided to illustrate the failure of $u^{(\epsilon)}\xrightarrow[\epsilon\rightarrow 0]{} u$ in $C([0,\infty);H^s(M))$. 

For $M=\mathbb{R}^d$, the linear evolution diverges. Let $u_0 = \langle \nabla\rangle^{-s}|\nabla|^{\frac{1-d}{2}}e^{-\frac{|x|^2}{2}}\in H^s(\mathbb{R}^d)$. An explicit computation yields
\begin{equation*}
\begin{split}
\| (U_\epsilon(t)-U(t))u_0 \|_{H^s}^2&=
    \int_{\mathbb{R}^d} 2(1-\cos(\epsilon^2t|\xi|^4))\langle \xi \rangle^{2s}|\widehat{u_0}|^2d\xi\\
    &=C_1\int_0^\infty (1-\cos(\epsilon^2 t r^4)) e^{-r^2}dr\\
    &=C_2\left(\frac{\sqrt{2 \pi } \left(\sin \left(\frac{1}{8} \left(\pi -\frac{1}{t \epsilon ^2}\right)\right) J_{\frac{1}{4}}\left(\frac{1}{8 t \epsilon ^2}\right)-\cos \left(\frac{1}{8} \left(\pi+\frac{1}{t \epsilon ^2} \right)\right) J_{-\frac{1}{4}}\left(\frac{1}{8 t \epsilon ^2}\right)\right)}{\sqrt{t} \epsilon }+4\right),
\end{split}
\end{equation*}
where $J_{\pm \frac{1}{4}}(x)$ is the Bessel function of the first kind and $C_1,C_2>0$ are dimensional constants. By direct computation, the last term is an increasing function in $t\in [0,\infty)$ whose limit as $t\rightarrow\infty$ is $4C_2$ for all $\epsilon>0$.

For $M=\mathbb{T}^d$, an exact solution to \eqref{adiabatic} with data $u_0(x)=\langle n \rangle^{-s}e^{in \cdot x}$ for $n\in \mathbb{Z}^d$ is given by  
\begin{equation*}
x\mapsto u^{(\epsilon)}(x,t):= e^{-it(\epsilon^2n^4+n^2-\langle n \rangle^{-2s})}\langle n \rangle^{-s}e^{inx}.    
\end{equation*}
Then,
\begin{equation*}
    \sup\limits_{t\in [0,\infty)} \| u^{(\epsilon)} - u \|_{H^s} = \sup\limits_{t\in [0,\infty)} c_d |e^{-it\epsilon^2 n^4}-1| = 2c_d,
\end{equation*}
for all $\epsilon>0$.
\end{remark}

\section{Acknowledgements.}

The author would like to appreciate his doctoral advisor Mark Kon for insightful comments.

\begin{appendices}
\section{Proof of \cref{mainprop}}
Consider the following cubic polynomial in $\xi_2$ with $\epsilon>0,\:\tau,\xi,\xi_1\in\mathbb{R}$:
\begin{equation*}
P(\xi_2)\coloneqq 4\epsilon^2\xi_2^3 + |\xi_1|^{2/3}(\epsilon^2\xi_1^2+2)\xi_2 + |(1+\epsilon^2(\xi_1-\xi)^2)(\xi_1-\xi)^2+\tau|.
\end{equation*}
Since $P$ has non-negative coefficients, it has a unique non-positive real root.
\begin{lemma}\label{technicallem}
Let $\xi>1, \tau\in (\frac{w(\xi)}{2},2w(\xi))$ and let $r(\xi_1)$ denote the unique negative root of $P$. Then,
\begin{enumerate}
    \item $r(\xi_1)=-\frac{|\xi_1|^{1/3}(\epsilon^2\xi_1^2+2)^{1/2}}{\sqrt{3}\epsilon}\sinh\bigg(\frac{1}{3}\sinh^{-1}\Big(3\sqrt{3}\epsilon \frac{(1+\epsilon^2(\xi_1-\xi)^2)(\xi_1-\xi)^2+\tau}{|\xi_1|(\epsilon^2\xi_1^2+2)^{3/2}}\Big)\bigg)$.
    \item $|r(\xi_1)|\gtrsim_\epsilon |\xi|^{\frac{4}{3}}$ for all $\xi_1 \in (-\infty,\infty)$ where the implicit constant does not depend on $\xi$.
    \item $|r(\xi_1)|\lesssim_\epsilon 
\begin{cases}
    w(\xi)^{1/3},& \text{if } |\xi_1|\in [0,\frac{\xi}{2}]\\
    |\xi_1|^{1/3}\langle \epsilon \xi_1 \rangle,              & \textit{otherwise}.
\end{cases}$
\end{enumerate}
\end{lemma}
\begin{proof}
The first statement is a hyperbolic trigonometric representation of a cubic root for a unique real root, which can be verified by direct substitution.

For the second statement, since $|r(-\xi_1)|\geq |r(\xi_1)|$ for all $\xi_1 \geq 0$ (due to $\xi>1$), it suffices to show $|r(\xi_1)|\gtrsim_\epsilon |\xi|^{4/3}$ for $\xi_1\geq 0$. Observe that $|r|$ is a decreasing function on $\xi_1\in (0,\xi)$ since $P(0)=\Big(1+\epsilon^2(\xi_1-\xi)^2\Big)(\xi_1-\xi)^2+\tau$ is decreasing and $\partial_{\xi_2}P(0) = |\xi_1|^{2/3}(\epsilon^2\xi_1^2+2)$ is increasing on $\xi_1\in(0,\xi)$; sketch a graph to see this. Hence for $\xi_1 \in (0,\xi]$,
\begin{equation*}
\begin{split}
|r(\xi_1)|&\gtrsim_{\epsilon} |\xi|^{1/3}\langle \epsilon \xi\rangle \sinh\bigg(\frac{1}{3}\sinh^{-1}\Big(3\sqrt{3}\epsilon\frac{\frac{w(\xi)}{2}}{|\xi|(\epsilon^2\xi^2+2)^{3/2}}\Big)\bigg)\gtrsim_\epsilon |\xi|^{4/3}.    
\end{split}
\end{equation*}
For $\xi_1 \in (\xi,2\xi]$,
\begin{equation}\label{phase}
\begin{split}
|r(\xi_1)|&\simeq_{\epsilon} |\xi_1|^{1/3}\langle \epsilon \xi_1\rangle \sinh\bigg(\frac{1}{3}\sinh^{-1}\Big(3\sqrt{3}\epsilon\frac{(1+\epsilon^2(\xi_1-\xi)^2)(\xi_1-\xi)^2+\tau}{|\xi_1|(\epsilon^2\xi_1^2+2)^{3/2}}\Big)\bigg)\\
    &\geq |\xi|^{1/3}\langle \epsilon \xi\rangle \sinh\bigg(\frac{1}{3}\sinh^{-1}\Big(3\sqrt{3}\epsilon\frac{\frac{w(\xi)}{2}}{2|\xi|(4\epsilon^2\xi^2+2)^{3/2}}\Big)\bigg)\\
    &\gtrsim_\epsilon |\xi|^{4/3},
\end{split}    
\end{equation}
where the last inequality of \eqref{phase} holds since for $\xi>1$,
\begin{equation*}
    \frac{\frac{w(\xi)}{2}}{2|\xi|(4\epsilon^2\xi^2+2)^{3/2}}\gtrsim \frac{w(\xi)}{\xi^4}\gtrsim_\epsilon 1.
\end{equation*}

Let $\xi_1 \in (2\xi,\infty)$. Since $\frac{\xi_1-\xi}{\xi_1}\geq \frac{1}{2}$, the argument inside $\sinh^{-1}$ is bounded below by a positive constant since
\begin{equation*}
    \frac{(1+\epsilon^2(\xi_1-\xi)^2)(\xi_1-\xi)^2+\tau}{|\xi_1|(\epsilon^2\xi_1^2+2)^{3/2}} \geq \frac{(1+\epsilon^2(\xi_1-\xi)^2)(\xi_1-\xi)^2}{|\xi_1|(\epsilon^2\xi_1^2+2)^{3/2}}\gtrsim_\epsilon \frac{(1+\epsilon^2(\xi_1-\xi)^2)(\xi_1-\xi)^2}{\xi_1^4}\gtrsim_\epsilon 1,
\end{equation*}
and this proves our claim as in \eqref{phase}.

To show the third statement, recall that $|r|$ is decreasing for $\xi_1\in [0,\xi)$, and therefore for such $\xi_1$
\begin{equation*}
    |r(\xi_1)|\leq |r(0)|\simeq_{\epsilon}(w(\xi)+\tau)^{1/3}\lesssim w(\xi)^{1/3}.
\end{equation*}
Furthermore for $\xi_1 \in [0,\frac{\xi}{2})$, we claim $|r(-\xi_1)|\lesssim_\epsilon |r(\xi_1)|$ from which $|r(\xi_1)|\lesssim_\epsilon w(\xi)^{1/3}$ follows for $|\xi_1|\in [0,\frac{\xi}{2})$. It suffices to show that there exists $C>0$ such that
\begin{equation}\label{phase2}
    \frac{\sinh\bigg(\frac{1}{3}\sinh^{-1}\Big(3\sqrt{3}\epsilon\frac{(1+\epsilon^2(\xi_1+\xi)^2)(\xi_1+\xi)^2+\tau}{|\xi_1|(\epsilon^2\xi_1^2+2)^{3/2}}\Big)\bigg)}{\sinh\bigg(\frac{1}{3}\sinh^{-1}\Big(3\sqrt{3}\epsilon\frac{(1+\epsilon^2(\xi_1-\xi)^2)(\xi_1-\xi)^2+\tau}{|\xi_1|(\epsilon^2\xi_1^2+2)^{3/2}}\Big)\bigg)}\leq C.
\end{equation}
Let $X,Y$ be the arguments inside the numerator and denominator of the hyperbolic sine of \eqref{phase2}, respectively. Then,
\begin{equation*}
    \frac{\sinh(X)}{\sinh(Y)} = \frac{e^X-e^{-X}}{e^Y-e^{-Y}}\leq \frac{e^X}{e^Y-e^{-Y}}\lesssim_\epsilon \frac{e^X}{e^Y},
\end{equation*}
where the last inequality follows since $Y$ is bounded below by a positive constant as can be observed from
\begin{equation*}
    \frac{(1+\epsilon^2(\xi_1-\xi)^2)(\xi_1-\xi)^2+\tau}{|\xi_1|(\epsilon^2\xi_1^2+2)^{3/2}}\geq \frac{\frac{w(\xi)}{2}}{\frac{|\xi|}{2}(\frac{\epsilon^2\xi^2}{4}+2)^{3/2}}\gtrsim_\epsilon 1.
\end{equation*}
Then using the identity
\begin{equation*}
    \sinh^{-1}(t) = \ln (t+\sqrt{1+t^2})
\end{equation*}
for all $t\in \mathbb{R}$ and letting $\alpha_1,\alpha_2$ be the arguments in the numerator and denominator of $\sinh^{-1}$ of \eqref{phase2}, respectively,
\begin{align*}
    \bigg(\frac{e^X}{e^Y}\bigg)^3&=\frac{\alpha_1+\sqrt{1+\alpha_1^2}}{\alpha_2+\sqrt{1+\alpha_2^2}}\leq \frac{\alpha_1+\sqrt{1+\alpha_1^2}}{\alpha_2}\lesssim_\epsilon \frac{\alpha_1}{\alpha_2},
\end{align*}
since $\alpha_1$ is bounded below by a positive constant (similar to $Y\gtrsim_{\epsilon} 1$). Note that our hypothesis on $\xi_1$ implies
\begin{equation*}
    -3\leq\frac{\xi_1+\xi}{\xi_1-\xi}\leq -1,
\end{equation*}
and therefore
\begin{equation*}
\begin{split}
\frac{\alpha_1}{\alpha_2}&=
    \frac{(1+\epsilon^2(\xi_1+\xi)^2)(\xi_1+\xi)^2+\tau}{(1+\epsilon^2(\xi_1-\xi)^2)(\xi_1-\xi)^2+\tau}\\
    &= \frac{(1+\epsilon^2(\xi_1+\xi)^2)(\xi_1+\xi)^2}{(1+\epsilon^2(\xi_1-\xi)^2)(\xi_1-\xi)^2+\tau}
    +\frac{\tau}{(1+\epsilon^2(\xi_1-\xi)^2)(\xi_1-\xi)^2+\tau}\\
    &\leq \frac{(1+\epsilon^2(\xi_1+\xi)^2)(\xi_1+\xi)^2}{(1+\epsilon^2(\xi_1-\xi)^2)(\xi_1-\xi)^2}+1\leq 9\Big(1+\Big(\frac{\xi_1+\xi}{\xi_1-\xi}\Big)^2\Big)+1 <\infty,
\end{split}    
\end{equation*}
as desired.

For $\xi_1>\frac{\xi}{2}$, we claim 
\begin{equation*}
\frac{(1+\epsilon^2(\xi_1-\xi)^2)(\xi_1-\xi)^2+\tau}{|\xi_1|(\epsilon^2\xi_1^2+2)^{3/2}}\lesssim\frac{(1+\epsilon^2(\xi_1-\xi)^2)(\xi_1-\xi)^2+w(\xi)}{|\xi_1|(\epsilon^2\xi_1^2+2)^{3/2}}\lesssim_\epsilon 1,    
\end{equation*}
which yields the claim immediately. Our hypothesis on $\xi_1$ implies
\begin{equation*}
    -1\leq \frac{\xi_1-\xi}{\xi_1}\leq 1,
\end{equation*}
from which
\begin{equation*}
\begin{split}
\frac{(1+\epsilon^2(\xi_1-\xi)^2)(\xi_1-\xi)^2+w(\xi)}{|\xi_1|(\epsilon^2\xi_1^2+2)^{3/2}}&=\frac{(1+\epsilon^2(\xi_1-\xi)^2)(\xi_1-\xi)^2}{|\xi_1|(\epsilon^2\xi_1^2+2)^{3/2}}+\frac{w(\xi)}{|\xi_1|(\epsilon^2\xi_1^2+2)^{3/2}}\\
    &\lesssim \frac{\Big(1+\epsilon^2(\xi_1-\xi)^2\Big)(\xi_1-\xi)^2}{\epsilon^3\xi_1^4}+c_\epsilon\lesssim \frac{1}{\epsilon^3\xi^2}+\frac{1}{\epsilon}+c_\epsilon\\
    &\lesssim_{\epsilon}1,    
\end{split}    
\end{equation*}
as desired. Arguing as above,
\begin{equation*}
|r(-\xi_1)|\lesssim_\epsilon |\xi_1|^{1/3}\langle\epsilon\xi_1\rangle    
\end{equation*}
for $\xi_1>\frac{\xi}{2}$ since for such $\xi_1$, $1\leq \frac{\xi_1+\xi}{\xi_1}\leq 3$.
\end{proof}
\begin{proof}[Proof of proposition \ref{mainprop}]
It suffices to prove the trilinear estimate in $X^{s,b}$.To illustrate this, fix $u\in X^{s,b}$ and let $\Tilde{u}\in X^{s,b}$ such that $u=\Tilde{u}$ on $t\in [0,T]$. Then since
\begin{equation*}
    \lVert N(u)\rVert_{X^{s+a,-\gamma}_{T}}\leq \lVert \eta(t/T)N(\Tilde{u})\rVert_{X^{s+a,-\gamma}}\lesssim_{\eta}\lVert N(\Tilde{u})\rVert_{X^{s+a,-\gamma}}\lesssim\lVert \Tilde{u} \rVert_{X^{s,b}}^3,
\end{equation*}
take the infimum over $\Tilde{u}$.

Let $u,v,w \in X^{s,b}$ and $\phi\in (X^{s+a,-\gamma})^*$. Define
\begin{equation*}
\begin{split}
f(\xi,\tau)&=|\hat{u}(\xi,\tau)|\langle \xi \rangle^s\langle \tau+w(\xi)\rangle^b,\:g(\xi,\tau)=|\hat{v}(\xi,\tau)|\langle \xi \rangle^s\langle \tau+w(\xi)\rangle^b,\\
    h(\xi,\tau)&=|\hat{w}(\xi,\tau)|\langle \xi \rangle^s\langle \tau+w(\xi)\rangle^b,\:\psi(\xi,\tau)=|\hat{\phi}(\xi,\tau)|\langle \xi \rangle^{-(s+a)}\langle \tau-w(\xi)\rangle^\gamma.
\end{split}
\end{equation*}
By the Plancherel's theorem, we have
\begin{equation*}
    |\langle (I-\epsilon^2\partial_{xx})^{-1}(u\overline{v})w,\phi\rangle_{L^2_{x,t}}|=|\langle \mathcal{F}[((I-\epsilon^2\partial_{xx})^{-1}(u\overline{v}))w],\mathcal{F}[\phi]\rangle_{L^2_{\xi, \tau}}|,
\end{equation*}
followed by the triangle inequality, which yields
\begin{equation*}
    \leq \int W\cdot f(\xi-\xi_1-\xi_2,\tau-\tau_1-\tau_2)g(-\xi_2,-\tau_2)h(\xi_1,\tau_1)\psi(\xi,\tau)d\xi_1d\xi_2d\xi d\tau_1 d\tau_2 d\tau,
\end{equation*}
where
\begin{equation*}
\begin{split}
W=&W(\xi,\xi_1,\xi_2,\tau,\tau_1,\tau_2)\\
    = &\frac{\langle\epsilon(\xi-\xi_1)\rangle^{-2}\langle \xi \rangle^{s+a}\langle \xi-\xi_1-\xi_2 \rangle^{-s}\langle \xi_2\rangle^{-s}\langle \xi_1 \rangle^{-s}}{\langle \tau-\tau_1-\tau_2+w(\xi-\xi_1-\xi_2)\rangle^b  \langle \tau_2-w(\xi_2)\rangle^b  \langle \tau_1 +w(\xi_1)\rangle^b  \langle \tau -w(\xi)\rangle^{\gamma}}.
\end{split}
\end{equation*}
We apply the Cauchy-Schwarz inequality (in variables $\xi_1,\xi_2,\tau_1,\tau_2$) on the product 
\begin{equation*}
    f(\xi-\xi_1-\xi_2,\tau-\tau_1-\tau_2)g(-\xi_2,-\tau_2)h(\xi_1,\tau_1)
\end{equation*}
and the rest of the integrand. The former, followed by the Young's inequality, yields
\begin{equation*}
    \lVert f^2\ast \overline{g^2}\ast h^2 \rVert_{L^1_{\xi, \tau}}^{1/2} \leq \lVert f \rVert_{L^2}\lVert g \rVert_{L^2}\lVert h \rVert_{L^2} = \lVert u \rVert_{X^{s,b}}\lVert v \rVert_{X^{s,b}}\lVert w \rVert_{X^{s,b}}.
\end{equation*}

The remaining part of the Cauchy-Schwarz inequality, followed by the $L^1-L^\infty$ H\"older's inequality yields
\begin{equation*}
    \lVert \psi \rVert_{L^2}\cdot\sup_{\xi,\tau}\Big(\int W^2 d\xi_1 d\xi_2 d\tau_1 d\tau_2\Big)^{\frac{1}{2}}.
\end{equation*}

Changing variable $\xi_1\mapsto \xi_1+\xi$, it suffices to show
\begin{equation}\label{finite}
    \sup_{\xi,\tau}\int \frac{\langle \epsilon\xi_1\rangle^{-4}\langle \xi \rangle^{2s+2a}\langle \tau -w(\xi)\rangle^{-2\gamma}\langle \xi_1+\xi_2 \rangle^{-2s}\langle \xi_2\rangle^{-2s}d\xi_1 d\xi_2 d\tau_1d\tau_2}{\langle \tau-\tau_1-\tau_2+w(\xi_1+\xi_2)\rangle^{2b}  \langle \tau_2-w(\xi_2)\rangle^{2b} \langle \xi_1+\xi \rangle^{2s} \langle \tau_1 +w(\xi_1+\xi)\rangle^{2b}}<\infty,
\end{equation}

where \eqref{finite} is reduced further by integrating in $\tau_1,\tau_2$:
\begin{equation}\label{finite2}
\begin{split}
\int \frac{ \langle \tau_1 +w(\xi_1+\xi)\rangle^{-2b}}{\langle \tau_1 -(\tau-\tau_2+w(\xi_1+\xi_2))\rangle^{2b}}d\tau_1&\lesssim \langle \tau_2-(\tau+w(\xi_1+\xi_2)+w(\xi_1+\xi))\rangle^{-2b},\\
    \int \frac{\langle \tau_2 - w(\xi_2)\rangle^{-2b}}{\langle \tau_2-(\tau+w(\xi_1+\xi_2)+w(\xi_1+\xi))\rangle^{2b}}d\tau_2&\lesssim \langle \tau+w(\xi_1+\xi_2)+w(\xi_1+\xi)-w(\xi_2) \rangle^{-2b}.
\end{split}
\end{equation}

Observing that $\langle \xi-A \rangle \langle \xi-B\rangle \gtrsim \langle A-B \rangle$ implies $\langle \xi_1+\xi_2\rangle^{2s}\langle \xi_2\rangle^{2s}\langle \xi_1+\xi \rangle^{2s}\gtrsim \langle \xi \rangle^{2s}$, \eqref{finite2} further reduces to showing
\begin{equation}\label{finite3}
    \sup\limits_{\xi,\tau}\frac{\langle\xi\rangle^{2a}}{\langle\tau-w(\xi)\rangle^{2\gamma}}\int \frac{\langle\epsilon \xi_1\rangle^{-4} d\xi_1 d\xi_2}{\langle\tau+w(\xi_1+\xi_2)+w(\xi_1+\xi)-w(\xi_2) \rangle^{2b}}<\infty.
\end{equation}

The expression $\tau+w(\xi_1+\xi_2)+w(\xi_1+\xi)-w(\xi_1) $ is a cubic polynomial in $\xi_2$ with an inflection point at $\xi_2 = -\frac{\xi_1}{2}$, and therefore after changing variable $\xi_2\mapsto \xi_2-\frac{\xi_1}{2}$, the integral of \eqref{finite3} becomes
\begin{equation*}
     \int \frac{\langle\epsilon \xi_1\rangle^{-4} d\xi_1 d\xi_2}{\langle4\epsilon^2\xi_1\xi_2^3 + \xi_1(\epsilon^2\xi_1^2+2)\xi_2 + (1+\epsilon^2(\xi_1+\xi)^2)(\xi_1+\xi)^2+\tau\rangle^{2b}}<\infty.
\end{equation*}

In doing the $\xi_2$-integral, if $\xi_1<0$, then via another change of variable $\xi_2\mapsto -\xi_2$, the integral is invariant when
\begin{equation*}
  4\epsilon^2\xi_1\xi_2^3 + \xi_1(\epsilon^2\xi_1^2+2)\xi_2 + (1+\epsilon^2(\xi_1+\xi)^2)(\xi_1+\xi)^2+\tau  
\end{equation*}

is replaced with
\begin{equation}\label{finite4}
    4\epsilon^2|\xi_1|\xi_2^3 + |\xi_1|(\epsilon^2\xi_1^2+2)\xi_2 + (1+\epsilon^2(\xi_1+\xi)^2)(\xi_1+\xi)^2+\tau.
\end{equation}

Similarly \eqref{finite4} can be replaced with
\begin{equation*}
     4\epsilon^2|\xi_1|\xi_2^3 + |\xi_1|(\epsilon^2\xi_1^2+2)\xi_2 + |(1+\epsilon^2(\xi_1+\xi)^2)(\xi_1+\xi)^2+\tau|,
\end{equation*} 
leaving the integral invariant.

Another change of variable $\xi_2\mapsto \frac{\xi_2}{|\xi_1|^{1/3}}$, followed by $\xi_1\mapsto -\xi_1$, eliminates the $\xi_1$-dependence in the leading coefficient of this cubic polynomial, and our task simplifies to showing the following expression is finite:
\begin{equation}\label{finite8}
    \sup\limits_{\xi,\tau}\frac{\langle\xi\rangle^{2a}}{\langle\tau-w(\xi)\rangle^{2\gamma}}
    \cdot\iint \frac{\langle\epsilon \xi_1\rangle^{-4}|\xi_1|^{-1/3} d\xi_1 d\xi_2}{\langle P(\xi_2) \rangle^{2b}}.
\end{equation}
It can be observed that $\sup\limits_{\xi\in\mathbb{R}}$ can be reduced to $\sup\limits_{\xi>0}$, which we assume henceforth; if $\xi<0$, let $\xi^\prime = -\xi$ and change variable $\xi_1\mapsto -\xi_1$ in the integral. However $\tau>0,\tau<0$ need to be discussed separately.

\textbf{Case I.} $\tau<0$.

Since $\langle \xi \rangle^{2a}\langle \tau-w(\xi)\rangle^{-2\gamma}\leq \langle \xi \rangle^{2a}\langle w(\xi)\rangle^{-2\gamma}\lesssim_\epsilon 1$, it suffices to show 
\begin{equation}\label{finite5}
    \sup_{\xi,\tau}\iint \frac{|\xi_1|^{-1/3}d\xi_1 d\xi_2}{\langle 4\epsilon^2\xi_2^3 + |\xi_1|^{2/3}(\epsilon^2\xi_1^2+2)\xi_2 + |(1+\epsilon^2(\xi_1-\xi)^2)(\xi_1-\xi)^2+\tau|\rangle^{2b}}<\infty.
\end{equation}
We do the $\xi_2$-integral in three regions: $(-\infty,r)\cup(r,0)\cup (0,\infty)$.\\

$(i)$. Consider the Taylor expansion of $|P(\xi_2)|$ on $(-\infty,r)$ at $\xi_2=r$.
\begin{equation*}
\begin{split}
|P(\xi_2)| &=-(12\epsilon^2r^2+|\xi_1|^{2/3}(\epsilon^2\xi_1^2+2))(\xi_2-r)-12\epsilon^2r(\xi_2-r)^2-4\epsilon^2(\xi_2-r)^3\\
    &\geq \max\Big(-|\xi_1|^{2/3}(\epsilon^2\xi_1^2+2)(\xi_2-r),-4\epsilon^2(\xi_2-r)^3\Big)\geq 0.
\end{split}
\end{equation*}
Integrate these lower bounds to obtain
\begin{equation*}
\begin{split}
|\xi_1|^{-1/3}\int_{-\infty}^r\frac{d\xi_2}{\langle P(\xi_2)\rangle^{2b}}&\leq|\xi_1|^{-1/3}\int_{-\infty}^r\frac{d\xi_2}{\langle |\xi_1|^{2/3}(\epsilon^2\xi_1^2+2)(\xi_2-r) \rangle^{2b}}\\
    &=|\xi_1|^{-1/3}\int_{0}^\infty\frac{d\xi_2}{\langle |\xi_1|^{2/3}(\epsilon^2\xi_1^2+2)\xi_2 \rangle^{2b}}\\
    &=\frac{|\xi_1|^{-1/3}}{|\xi_1|^{2/3}(\epsilon^2\xi_1^2+2)}\int_{0}^\infty\frac{d\xi_2}{\langle\xi_2 \rangle^{2b}}
    \simeq_b \frac{1}{|\xi_1|(\epsilon^2\xi_1^2+2)},
\end{split}
\end{equation*}
and
\begin{equation*}
    |\xi_1|^{-1/3}\int_{-\infty}^r\frac{d\xi_2}{\langle P(\xi_2)\rangle^{2b}}\leq|\xi_1|^{-1/3}\int_{-\infty}^r\frac{d\xi_2}{\langle 4\epsilon^2(\xi_2-r)^3\rangle^{2b}}=|\xi_1|^{-1/3}\int_{0}^\infty\frac{d\xi_2}{\langle 4\epsilon^2\xi_2^3\rangle^{2b}}\simeq_{b,\epsilon}|\xi_1|^{-1/3}.
\end{equation*}
Hence
\begin{equation*}
    |\xi_1|^{-1/3}\int_{-\infty}^r\frac{d\xi_2}{\langle P(\xi_2)\rangle^{2b}}\lesssim_{b,\epsilon}\min\Big(\frac{1}{|\xi_1|(\epsilon^2\xi_1^2+2)},|\xi_1|^{-1/3}\Big)
\end{equation*}
and \eqref{finite5} follows by integrating in $\xi_1$.

$(ii)$. Let $\xi_2 \in (r,0)$. Since $r$ is a root of $P(\xi_2)$,
\begin{equation*}
\begin{split}
|P(\xi_2)|&=P(\xi_2)=(\xi_2-r)(4\epsilon^2\xi_2^2+4\epsilon^2r\xi_2+|\xi_1|^{2/3}(\epsilon^2\xi_1^2+2)+4\epsilon^2r^2)\\
    &\geq (\xi_2-r)(|\xi_1|^{2/3}(\epsilon^2\xi_1^2+2))\geq 0.
\end{split}
\end{equation*}
Using this lower bound,
\begin{equation}\label{finite6}
\begin{split}
|\xi_1|^{-1/3}\int_r^0\frac{d\xi_2}{\langle P(\xi_2)\rangle^{2b}}&\leq |\xi_1|^{-1/3}\int_r^0 \frac{d\xi_2}{\langle (\xi_2-r)(|\xi_1|^{2/3}(\epsilon^2\xi_1^2+2))\rangle^{2b}}\\
    &=|\xi_1|^{-1/3}\int_0^{|r|}\frac{d\xi_2}{\langle \xi_2\cdot(|\xi_1|^{2/3}(\epsilon^2\xi_1^2+2))\rangle^{2b}}\\
    &\leq \frac{|\xi_1|^{-1/3}}{|\xi_1|^{2/3}(\epsilon^2\xi_1^2+2)}\int_0^\infty \frac{d\xi_2}{\langle \xi_2\rangle^{2b}}\simeq_{b}\frac{1}{|\xi_1|(\epsilon^2\xi_1^2+2)}.
\end{split}
\end{equation}
On the other hand,
\begin{equation*}
\begin{split}
P(\xi_2)&=(12\epsilon^2r^2+|\xi_1|^{2/3}(\epsilon^2\xi_1^2+2))(\xi_2-r)+12\epsilon^2r(\xi_2-r)^2+4\epsilon^2(\xi_2-r)^3\\
    &\geq 4\epsilon^2(\xi_2-r)^3\geq 0,
\end{split}
\end{equation*}
where the inequality holds since
\begin{equation*}
(12\epsilon^2r^2+|\xi_1|^{2/3}(\epsilon^2\xi_1^2+2))(\xi_2-r)+12\epsilon^2r(\xi_2-r)^2\geq 0,    
\end{equation*}
on $\xi_2 \in (r,0)$. Then,
\begin{equation}\label{finite7}
\begin{split}
|\xi_1|^{-1/3}\int_r^0\frac{d\xi_2}{\langle P(\xi_2)\rangle^{2b}}&\leq |\xi_1|^{-1/3}\int_r^0 \frac{d\xi_2}{\langle 4\epsilon^2(\xi_2-r)^3\rangle^{2b}}\\
    &=|\xi_1|^{-1/3}\int_0^{|r|}\frac{d\xi_2}{\langle 4\epsilon^2\xi_2^3\rangle^{2b}}\leq |\xi_1|^{-1/3}\int_0^{\infty}\frac{d\xi_2}{\langle 4\epsilon^2\xi_2^3\rangle^{2b}}\simeq_{b,\epsilon}|\xi_1|^{-1/3}.
\end{split}
\end{equation}
Hence
\begin{equation*}
    |\xi_1|^{-1/3}\int_{r}^0\frac{d\xi_2}{\langle P(\xi_2)\rangle^{2b}}\lesssim_{b,\epsilon}\min\Big(\frac{1}{|\xi_1|(\epsilon^2\xi_1^2+2)},|\xi_1|^{-1/3}\Big)
\end{equation*}
and the desired result follows by integrating with respect to $\xi_1$.

$(iii)$. Similarly on $\xi_2\in (0,\infty)$,
\begin{equation}
    |\xi_1|^{-1/3}\int_{0}^\infty\frac{d\xi_2}{\langle P(\xi_2)\rangle^{2b}}\lesssim_{b,\epsilon}\min\Big(\frac{1}{|\xi_1|(\epsilon^2\xi_1^2+2)},|\xi_1|^{-1/3}\Big),
\end{equation}
where the two lower bounds of $|P(\xi_2)|=P(\xi_2)\geq \max\Big(4\epsilon^2\xi_2^3,|\xi_1|^{2/3}(\epsilon^2\xi_1^2+2)\xi_2\Big)$ are used to argue as \eqref{finite6}, \eqref{finite7}. This concludes the proof for case I.

\textbf{Case II.} $\tau \in [0,\frac{w(\xi)}{2}]\cup [2w(\xi),\infty)$.

If $\tau \in [0,\frac{w(\xi)}{2}]$, then $w(\xi)-\tau \geq \frac{w(\xi)}{2}\geq 0$, and therefore $\langle \tau-w(\xi)\rangle^{-2\gamma}\leq \langle \frac{w(\xi)}{2}\rangle^{-2\gamma}$. If $\tau \in [2w(\xi),\infty)$, then $\tau-w(\xi)\geq w(\xi)$, and therefore $\langle \tau-w(\xi)\rangle^{-2\gamma}\leq \langle w(\xi)\rangle^{-2\gamma}$. In either case, the analysis reduces to showing
\begin{equation}
    \sup_{\xi,\tau}\iint \frac{|\xi_1|^{-1/3}d\xi_1 d\xi_2}{\langle 4\epsilon^2\xi_2^3 + |\xi_1|^{2/3}(\epsilon^2\xi_1^2+2)\xi_2 + (1+\epsilon^2(\xi_1-\xi))(\xi_1-\xi)^2+\tau\rangle^{2b}}<\infty,
\end{equation}
which can be done as in case I.

\textbf{Case III.} $\tau \in (\frac{w(\xi)}{2},2w(\xi))$. 

On this region, $\langle \tau - w(\xi)\rangle^{-2\gamma}\leq 1$. It is shown that the double integral of \eqref{finite8} is of $O(|\xi|^{-2a})$ as $|\xi|\rightarrow \infty$. As before, we derive lower bounds on $|P(\xi_2)|$ on three regions, $(-\infty,r)\cup (r,0)\cup (0,\infty)$. Moreover let $\xi>1$ by the Extreme Value Theorem since $\left\{(\xi,\tau):\xi \in [0,1], \tau \in [\frac{w(\xi)}{2},2w(\xi)]\right\}\subseteq \mathbb{R}^2$ is compact.

$(i)$. On $\xi_2\in (-\infty,r)$,
\begin{equation*}
\begin{split}
|P(\xi_2)| =-(12\epsilon^2r^2+|\xi_1|^{2/3}(\epsilon^2\xi_1^2+2))&(\xi_2-r)-12\epsilon^2r(\xi_2-r)^2-4\epsilon^2(\xi_2-r)^3\\
    &\geq -(12\epsilon^2r^2+|\xi_1|^{2/3}(\epsilon^2\xi_1^2+2))(\xi_2-r)\geq 0,
 \end{split}
\end{equation*}
from which
\begin{equation*}
\begin{split}
|\xi_1|^{-1/3}\int_{-\infty}^r\frac{d\xi_2}{\langle P(\xi_2)\rangle^{2b}}&\leq |\xi_1|^{-1/3}\int_0^\infty \frac{d\xi_2}{\langle( 12\epsilon^2r^2+|\xi_1|^{2/3}(\epsilon^2\xi_1^2+2))\xi_2\rangle^{2b}}\\
    &\simeq_{b}\frac{|\xi_1|^{-1/3}}{12\epsilon^2r^2+|\xi_1|^{2/3}(\epsilon^2\xi_1^2+2)}.
\end{split}
\end{equation*}
Change variable $z=\xi_1^{\frac{8}{3}}+c_\epsilon \xi^{\frac{8}{3}}$ where $c_\epsilon>0$ is to be determined and $p\in [1,3)$. Then,
\begin{equation}\label{finite9}
\begin{split}
&\int_0^\infty \frac{\langle \epsilon\xi_1\rangle^{-4}|\xi_1|^{-1/3}d\xi_1}{12\epsilon^2r^2+|\xi_1|^{2/3}(\epsilon^2\xi_1^2+2)}\\
    =&
    \int_0^\xi \frac{\langle \epsilon\xi_1\rangle^{-4}|\xi_1|^{-1/3}d\xi_1}{12\epsilon^2r^2+|\xi_1|^{2/3}(\epsilon^2\xi_1^2+2)}+\int_\xi^\infty \frac{\langle \epsilon\xi_1\rangle^{-4}|\xi_1|^{-1/3}d\xi_1}{12\epsilon^2r^2+|\xi_1|^{2/3}(\epsilon^2\xi_1^2+2)}\\
    \lesssim_\epsilon & |\xi|^{-\frac{8}{3}}\lVert \langle\epsilon \xi_1\rangle^{-4}\rVert_{L^{p^\prime}_{(0,\xi)}}\lVert \xi_1^{-1/3}\rVert_{L^p_{(0,\xi)}}+\langle \epsilon\xi\rangle^{-4}\xi^{-1/3}\int_\xi^\infty \frac{d\xi_1}{\xi_1^{8/3}+c_\epsilon \xi^{8/3}}\\
    \lesssim_p & \langle\xi\rangle^{-(3-\frac{1}{p})}+\langle \epsilon\xi\rangle^{-4}\xi^{-1/3}\int_{c_\epsilon\xi^{8/3}}^\infty \frac{dz}{z(z-c_\epsilon\xi^{8/3})^{\frac{5}{8}}}\lesssim \langle \xi \rangle^{-(3-\frac{1}{p})},
\end{split}    
\end{equation}
where the lower bound of $|r|$ is by \cref{technicallem}.

$(ii)$. On $\xi_2\in (r,0)$, 
\begin{equation*}
\begin{split}
|P(\xi_2)|&=P(\xi_2)\\
&\geq \frac{(1+\epsilon^2(\xi_1-\xi)^2)(\xi_1-\xi)^2+\tau}{|r|}\xi_2+(1+\epsilon^2(\xi_1-\xi)^2)(\xi_1-\xi)^2+\tau\geq 0,
\end{split}
\end{equation*}
and therefore
\begin{equation*}
    \int_r^0 \frac{|\xi_1|^{-1/3} d\xi_2}{\langle \frac{(1+\epsilon^2(\xi_1-\xi)^2)(\xi_1-\xi)^2+\tau}{|r|}\xi_2+(1+\epsilon^2(\xi_1-\xi)^2)(\xi_1-\xi)^2+\tau\rangle^{2b}}\lesssim_b \frac{|\xi_1|^{-1/3}|r(\xi_1)|}{\langle (1+\epsilon^2(\xi_1-\xi)^2)(\xi_1-\xi)^2 + \frac{w(\xi)}{2}\rangle}.
\end{equation*}
We change variable $z=\xi_1^4 + \frac{\xi^4}{2}$ and integrate in $\xi_1$ as in case $(i)$, and use the upper bound of \cref{technicallem} to obtain
\begin{equation}\label{finite10}
\begin{split}
&\int_{-\infty}^\infty\int_r^0\frac{\langle \epsilon\xi_1\rangle^{-4}|\xi_1|^{-1/3}d\xi_2 d\xi_1}{\langle P(\xi_2)\rangle^{2b}}\\
    \lesssim& \int_0^{\xi}\frac{\langle \epsilon\xi_1\rangle^{-4}|\xi_1|^{-1/3}w(\xi)^{1/3}d\xi_1}{\langle\frac{w(\xi)}{2}\rangle}+\int_{\xi}^\infty \frac{\langle \epsilon\xi_1\rangle^{-3} d\xi_1}{\langle (1+\epsilon^2(\xi_1-\xi)^2)(\xi_1-\xi)^2+\frac{w(\xi)}{2}\rangle}\\
    \lesssim_b &\langle \xi \rangle^{-(3-\frac{1}{p})}+\frac{\langle\epsilon\xi\rangle^{-3}}{\epsilon^2}\int_0^\infty \frac{ d\xi_1}{\xi_1^4+\frac{\xi^4}{2}}\lesssim\langle \xi \rangle^{-(3-\frac{1}{p})}.    
\end{split}
\end{equation}
$(iii)$. On $\xi_2\in (0,\infty)$,
\begin{equation*}
    |P(\xi_2)|=P(\xi_2)\geq 4\epsilon^2\xi_2^3+\frac{w(\xi)}{2}\geq 0.
\end{equation*}
Change variable $z=4\epsilon^2\xi_2^3+\frac{w(\xi)}{2}$ to obtain
\begin{equation*}
\begin{split}
|\xi_1|^{-1/3}\int_0^\infty\frac{d\xi_2}{\langle P(\xi_2)\rangle^{2b}}&\leq |\xi_1|^{-1/3}\int_0^\infty\frac{d\xi_2}{\langle 4\epsilon^2\xi_2^3+\frac{w(\xi)}{2}\rangle^{2b}}\\
    &\simeq_{\epsilon} |\xi_1|^{-1/3}\int_{\frac{w(\xi)}{2}}^\infty \frac{dz}{\langle z\rangle^{2b}(z-\frac{w(\xi)}{2})^{2/3}}\lesssim_b |\xi_1|^{-\frac{1}{3}}|w(\xi)|^{-\frac{2}{3}},
\end{split}
\end{equation*}
where in the last inequality, we note that $\langle z\rangle^{-2b}\lesssim |z|^{-1}$ on the region of integration. On the other hand, we use $4\epsilon^2\xi_2^3+(1+\epsilon^2(\xi_1-\xi)^2)(\xi_1-\xi)^2+\frac{w(\xi)}{2}\geq 0$ as another lower bound to derive a similar estimate for $|\xi_1|>2\xi$ on which $|\xi_1-\xi|\geq \frac{|\xi_1|}{2}$.
\begin{equation*}
\begin{split}
|\xi_1|^{-1/3}\int_0^\infty\frac{d\xi_2}{\langle P(\xi_2)\rangle^{2b}}&\lesssim_{b,\epsilon}\frac{|\xi_1|^{-1/3}}{\langle (1+\epsilon^2(\xi_1-\xi)^2)(\xi_1-\xi)^2+\frac{w(\xi)}{2}\rangle^{2/3}}\leq \frac{|\xi_1|^{-1/3}}{\langle\frac{\epsilon^2\xi_1^4}{16}+\frac{w(\xi)}{2}\rangle^{2/3}}\\
    &\lesssim \frac{|\xi|^{-1/3}}{\Big(\frac{\epsilon^2\xi_1^4}{16}+\frac{w(\xi)}{2}\Big)^{2/3}},    
\end{split}
\end{equation*}
and the remaining $\xi_1$ integral proceeds as \eqref{finite9}, \eqref{finite10}.
\end{proof}
\end{appendices}

\bibliographystyle{abbrv}
\bibliography{ref3}
\end{document}